\documentclass{cimart}

\usepackage{amsfonts,amssymb,amscd}
\usepackage{amssymb}
\usepackage{graphics}
\usepackage{mathrsfs}
\usepackage{graphicx}
\usepackage{graphics}
\usepackage{textcomp}
\usepackage{gensymb}

\usepackage{epsf}
\usepackage{epsfig}
\usepackage{graphicx}

\usepackage{amsmath}
\usepackage{amssymb}
\usepackage{amsfonts}
\usepackage{pifont}
\usepackage{gensymb}

\usepackage{xcolor}

\newcommand{\p}{\mathbb{P}}

\newcommand{\dl}{\lambda}
\newcommand{\da}{\alpha}
\newcommand{\de}{\varepsilon}
\newcommand{\db}{\beta}

\newcommand{\C}{\mathbb{C}}
\newcommand{\R}{\mathbb{R}}
\newcommand{\Z}{\mathbb{Z}}
\newcommand{\N}{\mathbb{N}}

\newcommand{\Fcal}{\mathcal{F}}
\newcommand{\calO}{{\mathcal{O}}}

\newcommand{\mcf}{\mathcal{D}}

\newcommand{\tf}{\widetilde{\mathcal{D}}}
\newcommand{\tXX}{{\widetilde{X}}}

\newcommand{\tilf}{\widetilde{\mathcal{F}}}
\newcommand{\fol}{\mathcal{F}}

\newcommand{\calp}{{\mathcal{P}}}
\newcommand{\calh}{\mathcal{H}}

\newcommand{\diff}{{{\rm Diff}\, ({\mathbb C}, 0)}}
\newcommand{\diffd}{{{\rm Diff}\, ({\mathbb C}^2, 0)}}
\newcommand{\diffdd}{{{\rm Diff}_1\, ({\mathbb C}^2, 0)}}
\newcommand{\diffn}{{{\rm Diff}\, ({\mathbb C^n}, 0)}}

\newcommand{\diffalpha}{{{\rm Diff}_{\alpha} ({\mathbb C}, 0)}}

\newcommand{\x}{{\boldsymbol x }}

\newtheorem{defi}{Definition}
\newtheorem{teo}{Theorem}
\newtheorem{prop}{Proposition}

\newtheorem{rem}{Remark}
\newtheorem{quest}{Question}

\title{Complex ODE{\small s}, singularity theory and dynamics}

\author{Helena Reis}

\authorinfo[
Helena Reis]{
	Mathematics Center of the University of Porto,
	Faculty of Economics of the University of Porto,
	Portugal}{%
	hreis@fep.up.pt
}

\abstract{%
These notes are a slightly enlarged version of my habilitation thesis, where our research interests and main results in the past few years are summarized. Most of the discussion revolves around complex ordinary differential equations and their underlying foliations, singularity theory and dynamical systems. Compared to the original text, a section containing some background material on holomorphic foliations was added. Also, some new results obtained in the past three years that are in line with the one presented in the habilitation were included.
}



\VOLUME{32}
\YEAR{2024}
\ISSUE{3}
\firstpage{267}
\DOI{https://doi.org/10.46298/cm.12651}

\begin{document}

\section{Introduction}

The origin of these notes goes back to my ``Habilitation thesis'' presented at the University of Porto in 2021, where we are supposed to present a description of an important component of our research in the past few years. In this sense, these notes have a very significant overlap with the text of my ``Habilitation thesis'' (available from~\cite{Reis_habilitation}), although the two texts do not coincide. Most notably, compared to~\cite{Reis_habilitation}, a section with some background material in the area has been added as well as some new results obtained in the past three years. Some additional minor modifications were also made. It should also be mentioned that Sections~\ref{sec:invariant_sets}, \ref{sec:localaspects} and~\ref{Sec:resolution} have a non-trivial intersection with the discussion conducted in~\cite{Handbook}.

Roughly speaking my research concerns the dynamics and the geometry of complex ordinary differential equations. More precisely, a good part of my research has been focused on local and global aspects of holomorphic vector fields and/or foliations on $3$-dimensional
(regular) manifolds. It is well known that studying the singularities of holomorphic foliations in dimension at least~$3$ is much harder than the analogous problem in dimension~$2$. Let us then begin by singling out some global difficulties arising in these problems that have no $2$-dimensional counterpart. In what follows, unless otherwise mentioned, (singular) holomorphic foliations are always of dimension~$1$. In other words, these foliations
are locally given by the (local) orbits of a holomorphic vector field having a singular set of codimension at least~$2$.

One of the main difficulties in the study of singularities of holomorphic foliations on ambient space of dimension at least~$3$ comes from the fact that these singularities encode some {\it global dynamics}\, on the divisors naturally associated with them. To explain the role played by these dynamics, we may think of the one-point
blow-up of a ``generic'' homogeneous vector field on $\C^3$. While the singularities of the blown-up foliation become ``simple'', the
understanding of the initial singularity clearly requires the understanding of the global foliation induced on the projective space
identified with the exceptional divisor by the homogeneous vector field in question. In general, this foliation possesses very
complicated dynamics leaving no algebraic ``object'' invariant. This phenomenon does not occur in dimension~$2$ since the foliation
induced on the projective line consists of the union of a unique leaf with finitely many singular points. Thus the dynamics obtained
on the divisor is rather trivial.

Another well-known additional difficulty in problems involving singularities in dimension greater than~$2$ is the absence of a
desingularization procedure as effective as Seidenberg's theorem valid in dimension~$2$. In fact, according to Seidenberg, for
every holomorphic foliation on a complex surface, there exists a finite sequence of one-point blow-ups such that the corresponding
transform of the initial foliation possesses only {\it elementary singular points}. Recall that a singular point is said to be {\it elementary}\,
if the foliation admits at least one eigenvalue different from zero at it. It turns out, however, that a faithful analogue of Seidenberg's
result for foliations on 3-manifolds cannot exist: there are some non-simple singularities that are persistent under blow-up transformations
(cf. \cite{C-R-S} as well as Section~\ref{Sec:resolution} of the present paper for details). Nonetheless different sorts of {\it final models}\, for certain ``desingularization'' procedures were
described for example in the following papers \cite{C-R-S, MQ-P, P}, and more recently, in \cite{RR_Resolution}.

We can also include in this ``list of additional difficulties'' some problems related to divergent normal forms (irregular singularities).
In the case of saddle-node singularities in $n$-dimensional manifolds, where $n \geq 3$, the foliation may admit two or more eigenvalues
equal to zero. In this case not only their formal normal forms are poorly understood but also the resummation techniques are much less
developed. Naturally, the same problem occurs for every other singularity whose rank of resonance relations is at least~$2$.

My research presented in this text is contained in the papers \cite{MRR, P-R, RR_separatrix, RR_stabilizers, RR_applications, RR_secondjet, RR, RR_Resolution, Reis06, Reis08}. Many of these papers
solve long-standing problems in the area, including the existence of invariant analytic surfaces for commuting vector fields, the
topological type of leaves associated with Arnold's $A^{2n+1}$ singularities, the proof that complete integrability is not a topological
invariant of 1-dimensional foliations on $(\C^3,0)$, a topological characterization of virtually solvable subgroups of ${\rm Diff} \,
(\C^2,0)$ and, more recently, the proof of a sharp resolution theorem for singularities of complete vector fields in dimension 3.

My research, however, also includes a distinguished class of singularities of vector fields, namely the \emph{semicomplete} (singularities of ) vector fields. Roughly speaking, a vector field is semicomplete if the solutions of the differential equation associated with it are all univalued. Furthermore, semicomplete singularities are the only ones that can be realized by a complete vector field on some complex manifold. Understanding the above-mentioned class of vector fields, both at the global level and at the level of germs, is a problem with interesting applications. As an example of application, we will see in Section~\ref{sec:localaspects} that results on singularities of semicomplete vector fields yield insight into some problems about bounds for the dimension of the automorphism group of compact complex manifolds. Another motivation to study these vector fields and their singular
points stems from the very fact that the semicomplete property is somehow akin to the Painlev\'e property for differential equations, albeit the two notions are not equivalent. As a matter of fact, as it happens with Painlev\'e property, semicomplete vector fields are also largely present - sometimes implicitly - in the literature of Mathematical Physics.

It should be noted that the semicompleteness is not an intrinsic property of foliations in the sense that we may have two vector fields inducing the same foliation, with one of them being semicomplete while the other is not so. This is an extra difficulty compared to the simple study of the foliations induced by the vector fields.

These notes are organized so that all background necessary for a certain section was already discussed/introduced in a previous one. The paper is then structured as follows. To begin, in section~\ref{sec:basics} we introduce the standard terminology along with the basic tools in the study of holomorphic foliations. Next, in Section~\ref{sec:invariant_sets} we will discuss the fundamental problem of existence of separatrices, i.e. existence of germs of analytic sets that are invariant by (germs of) singular foliations and passing through a singular point. The existence of separatrices for foliations on $(\C^2,0)$ has been established by Camacho and Sad in 1982 (cf.~\cite{CS}) and counterexamples to their existence were provided some years later for $1$-dimensional foliations, while counterexamples in the codimension-$1$ case was established earlier (cf. \cite{GM-L} and \cite{Joa}, respectively). In this section we discuss the counterexamples in question and present a way to construct many other counterexamples in the case of codimension~$1$ foliations. The idea of this construction passes through the fact that a ``generic'' foliation on $\C\p(2)$ has no invariant algebraic curve. We then explain how we were able to establish the existence of separatrices for codimension~$1$ foliation on $(\C^3,0)$ that are spanned by two commuting vector fields by exploiting the dynamics of the associated foliation on the corresponding exceptional divisor.

Sections~\ref{sec:sc_global_dynamics} and~\ref{sec:localaspects} are devoted to the special class of semicomplete vector fields and their corresponding singularities. Although this class is quite ``small'' in an appropriate sense, its importance in terms of properties and applications completely justifies its study, as it will be made clear in the mentioned sections. Each one of these two sections have different nature. In the first one, we discuss global aspects associated with the mentioned vector field (or, more precisely, with a suitable subclass of it - namely, the class of complete vector fields on $(\C^n,0)$), while in the second section we focus on their local aspects and applications to the problem of bounds for the automorphism group of a compact complex manifold.

Section~\ref{Sec_topology} is the only section concerning exclusively foliations in a two-dimensional space. In this section we present the geometric study of the foliations associated to Arnold singularities $A^{2n+1}$ along with the results on pseudogroups of ${\rm Diff} \, (\C,0)$ obtained in~\cite{MRR} and~\cite{RR_stabilizers} that allowed us to establish a long standing question about the topological type of leaves associated with Arnold singularities.

Still in the context of pseudogroups on ${\rm Diff} \, (\C,0)$, a characterization of those having only finite orbits was presented by Mattei and Mattei-Moussu. As shown by Mattei and Moussu, such pseudogroups are strictly related with integrability of the groups and of the foliations having such groups as holonomy group. Extensions of Mattei and Moussu results are discussed on Section~\ref{sec_integrability}.

Finally, Section~\ref{Sec:resolution} is devoted to resolution theorems for $1$-dimensional foliations on $(\C^3,0)$ (recall that final models for resolution procedures of codimension-$1$ foliations on $(\C^3,0)$ are well understood after \cite{cano}). We conduct a thorough discussion about the content of the resolution theorems deduced in~\cite{MQ-P} and in~\cite{RR_Resolution}, highlighting virtues and potential limitations.


\section{Basics in the local theory of holomorphic foliations}\label{sec:basics}

Let us then begin this section by recalling the definition of (singular) holomorphic foliations.

\begin{defi}\label{defi:foliation}
	Let $M$ be a complex manifold of dimension~$m$. A singular holomorphic foliation of dimension~$n$ consists of a distinguished coordinate covering $\fol = \{(U_i, \varphi_i)\}$, where $\varphi_i: U_i \to \C^m$, satisfying the conditions below
	\begin{itemize}
		\item[1.] $U = \{U_i\}$ is an open cover of $M \setminus S$, where $S$ is an analytic subset of $M$ of codimension at least~$2$ (possibly empty);
		
		\item[2.] whenever $U_i \cap U_j \ne \emptyset$, the diffeomorphism $\varphi_{ij}$ defined as $\varphi_{ij} = \varphi_i \circ \varphi_j^{-1}: \varphi_j (U_i \cap U_j) \to \varphi_i (U_i \cap U_j)$ takes on the form
		\[
		\varphi_{ij} (x,y) = (f_{ij}(x,y), g_{ij}(y)) \, ,
		\]
		with $x \in \C^{n}$ and $y \in \C^{m-n}$.
	\end{itemize}
\end{defi}

For simplicity, and if no misunderstanding is possible, a foliation  $\fol = \{(U_i, \varphi_i)\}$ will be simply denoted by $\fol$. The set $S$ in Definition~\ref{defi:foliation} is said to be the \emph{singular set} of the foliation,.

Let $\fol = \{(U_i, \varphi_i)\}$ be a foliation on a complex manifold $M$. A {\it plaque} of $\fol$ is a subset of $M$ given as $\varphi_i^{-1}(\{y = {\rm cte}\})$, for some $i$ and some constant ${\rm cte}$. The plaques of $\fol$ define a relation $\sim$ on $M$ as follows: for every $x, \, y \in M$ we say that $x \sim y$ if and only if there exists a finite sequence of plaques $\alpha_1, \, \ldots, \alpha_k$ such that
\begin{itemize}
	\item $x \in \alpha_1$;
	\item $y \in \alpha_k$;
	\item $\alpha_i \cap \alpha_{i+1} \ne \emptyset$, for every $1 \leq i \leq k-1$.
\end{itemize}
A {\it leaf} of the foliation $\fol$ is an equivalence class of the relation $\sim$. Furthermore, the leaf (or, equivalently, the equivalence class) of $\fol$ thought a point $p \in M$ will be denote by $L_p$.

Since we will focous on foliations of dimension~$1$, it should be recalled that such foliations may equivalently be defined by means of holomorphic vector fields. In fact, we have the following

\begin{defi}\label{defi:foliations_dim1}
	A singular holomorphic foliation of dimension~$1$ defined on a complex manifold $M$ consists of the following data:
	\begin{itemize}
		\item[1.] an atlas ${(U_i, \varphi_i)}$ compatible with the complex structure on $M$, with $\varphi_i: U_i \to V_i \subseteq \C$;
		\item[2.] a holomorphic vector field $X_i$ defined on $V_i$, for each $i$, with singular set of codimension at least~$2$;
		\item[3.] if $U_i \cap U_j \ne \emptyset$, then
		\[
		(\varphi_j \circ \varphi_i)^{-1}_{\ast} X_i (\varphi_i(U_i \cap U_j)) = h_{ij} (x_1, \ldots, x_n) \, X_j (\varphi_j (U_i \cap U_j))
		\]
		for some no-where holomorphic function $h_{ij}: \varphi_i(U_i \cap U_j) \to \C$.
	\end{itemize}
\end{defi}

Note that, in Definition~\ref{defi:foliations_dim1}, $X_i$ and $X_j$ need not to coincide on $U_i \cap U_j$. They just need to induce the same direction at every single point of $U_i \cap U_j$. In the case where $h_{ij}$ is constant equal to~$1$ for every $(i,j)$, then $\{X_i\}$ defines an actual vector field on $M$.

With respect to the present definition, the singular set ${\rm Sing} \, \fol$ of $\fol$ is then defined as the union over $i$ of the sets $\varphi_i^{-1} ({\rm Sing} \, (X_i))$ on $M$. Thus, there immediately follows that the singular set of any holomorphic $1$-dimensional foliation has codimension at least two and that a foliation has no divisor of zeros.

If we are given a holomorphic vector field $X$ on a complex manifold, with singular set of codimension at least~$2$, the leaves of the foliation induced by $X$ is nothing but the integral curves of the differential equation associated with it. Note, however, that the integral curves of a vector field taking on the form $Y = fX$, for some holomorphic function $f$, coincide with the integral curves of $X$ away from the zero divisor of $Y$, i.e. away from $\{f=0\}$. In this case, we say that $X$ and $Y$ induce the foliation. A vector field inducing a foliation $\fol$ and having singular set of codimension at least~$2$ is said to be a {\it representative} of $\fol$.

In turn, a codimesion~$1$ foliation can be defined by means of an integrable $1$-form, i.e. a $1$-form $\omega$ such that $\omega \wedge d\omega \equiv 0$. Recall that the kernel of  a (holomorphic) $1$-form $\omega$ defines, away from its singular set, a distribution of complex hyperplanes and this distribution is integrable if and only if $\omega \wedge d\omega \equiv 0$. Summarizing.

\begin{defi}
	A singular holomorphic foliation of codimension~$1$ on a complex manifold $M$ consists of the following data:
	\begin{itemize}
		\item[1.] an atlas ${(U_i, \varphi_i)}$ compatible with the complex structure on $M$, with $\varphi_i: U_i \to V_i \subseteq \C$;
		\item[2.] a collection of differential $1$-forms $\omega_i$ defined on $V_i$, for each $i$, with singular set of codimension at least~$2$ and such that $\omega_i \wedge d\omega_i$ vanishes identically;
		\item[3.] if $U_i \cap U_j \ne \emptyset$, then
		\[
		(\varphi_j \circ \varphi_i)^{-1}_{\ast} \omega_i = h_{ij} \, \omega_j
		\]
		for some no-where holomorphic function $h_{ij}: \varphi_i(U_i \cap U_j) \to \C$.
	\end{itemize}
\end{defi}

The notion of {\it representative $1$-form} for a codimension-$1$ foliation can be defined analogously to the notion of representative vector field for $1$-dimensional foliations.

One of the basic objects in the study of local theory of foliations is the so-called {\it separatrix}.

\begin{defi}\label{defi:separatrix}
	Let $\fol$ be a foliation of dimension~$k$ on $(\C^n, 0)$. A separatrix for $\fol$ is the germ of an irreducible analytic set $S$ of dimension~$k$ passing through the origin and invariant by $\fol$.
\end{defi}

Separatrices are objects of natural interest since they fit the framework of ``invariant manifolds'' in dynamical systems and their presence yields specific solutions for the vector field in question that can be understood in detail. Section~\ref{sec:invariant_sets} is entirely devoted to results on existence of separatrices.

Another basic tool that is particularly useful to understand the behavior of foliations at singular points is the blow-up. Let us begin by recalling the notion of blow-up centred at a point. We begin with the case of the affine space $\C^n$.

\begin{defi}\label{def:blowup}
	The blow-up of $\C^n$ at the origin is a complex manifold $\widetilde{\C}^n$ obtained by identifying $n$ copies of $\C^n$ in the following way. If $(x_1, \ldots, x_n)$ stands for local coordinates on $\C^n$, the $n$ charts on $\widetilde{\C}^n$ can be defined as $(u_1, \ldots, u_{i-1},v_i, u_{i+1}, \ldots, u_n)$ through the relations
	\begin{align*}
		\begin{cases}
			x_i &= v_i \\
			x_j &= v_i \, u_j \, , \quad \text{for all} \, j \ne i \, .
		\end{cases}
	\end{align*}
\end{defi}

The blow-up mapping $\pi: \widetilde{\C}^n \to \C^n$ is given on the different charts above by \[\pi(u_1, \ldots, u_{i-1},v_i, u_{i+1}, \ldots, u_n) = (u_1 v_i, \ldots, u_{i-1} v_i,v_i, u_{i+1} v_{i+1}, \ldots, u_n v_i).\] Moreover, it verifies the following:
\begin{itemize}
	\item $\pi^{-1} (0)$ is a well defined submanifold of $\widetilde{\C}^n$, isomorphic to the projective space $\C \p(n-1)$, that is locally given by $\{v_i = 0\}$;
	
	\item the restriction of $\pi$ to $\widetilde{\C}^n \setminus E$ ($\pi: \widetilde{\C}^n \setminus E \to \C^n \setminus 0$) is a holomorphic diffeomorphism;
	
	\item $\pi$ is proper, i.e. the pre-image of a compact set is also compact.
\end{itemize}

The definition of blow-up of a complex manifold $M$ centered at a point $p \in M$ can be obtained by means of the preceding construction, by using the complex coordinates of the manifold in question. To be more precise, consider a complex manifold $M$ and fix a point $p \in M$. Consider a local coordinate chart $\psi: U \to W \subseteq \C^n$ defined on a neighborhood $U$ of $p$ and such that $\psi (p) = 0$.
Let $\widetilde{W}$ stands for the preimage of $W$ through $\pi$, where $\pi$ stands for the blow-up map of $\C^n$ centered at the origin. Let then $M'$ be the disjoint union of $M \setminus \widetilde{W}$, and consider the following equivalence relation. Fix points $p_0 \in U \setminus \{p\}$ and $p_1 \in \widetilde{W} \setminus E$. We have
\[
q_0 \sim q_1 \quad \Leftrightarrow \quad q_1 = \pi^{-1} (\psi (q_0)) \, .
\]
The blow-up $\widetilde{M}$ of $M$ at $p$ is defined as the quotient of $M'$ by this equivalence relation, $M'/\sim$. Note that $\widetilde{M}$ is indeed a smooth complex manifold for $\widetilde{W}$ is a manifold and $\pi^{-1} \circ \psi \, : U \setminus \{p\} \to \widetilde{W} \setminus E$ is a holomorphic
diffeomorphism.
Similarly, there is a blow-up mapping from $\widetilde{M}$ to $M$ (which will also be denoted by $\pi$) that is proper and takes $E$ to $p$, i.e. $\pi(E) = p$. Moreover, the restriction of $\pi$ to $\widetilde{M} \setminus E$, $\pi : \widetilde{M} \setminus E \to M \setminus \{p\}$ is a holomorphic diffeomorphism.

Blow-ups centered at higher dimensional submanifolds of $M$ can also be defined. However, since the use we made of them is limited, we content ourselves to refer to~\cite{Shaf} for accurate definitions.

Given a foliation $\fol$, the transform of $\fol$ under a blow-up map is called the {\it blow-up of $\fol$} and it will usually be denoted by $\tilf$. Recall, however, that the blow-up space contains an exceptional divisor and the latter may or may not be invariant by the transformed foliation. This issue gives rise to the notion of {\it dicritical foliation}.

\begin{defi}
	Let $\fol$ be a holomorphic foliation on a complex manifold $M$. Consider a blow-up map $\pi : \widetilde{M} \to M$ centered at a subset $C$, $C \subseteq {\rm Sing} \, (\fol)$, where ${\rm Sing} \, (\fol)$ stands for the singular set of $\fol$. The foliation $\fol$ is said to be dicritical with respect to $\pi$ if the transformed foliation $\tilf$ does not leave the exceptional divisor $\pi^{-1} \, (C)$ invariant.
\end{defi}


Recall that on the study of holomorphic foliations, it is often necessary to iterate blow-ups. The notion of {\it dicritical foliation} can thus be made intrinsic as follows.

\begin{defi}
	A holomorphic foliation $\fol$ is called dicritical if there is a finite sequence of blow-ups with invariant centers such that the total exceptional divisor possesses an irreducible component that is not invariant for the transform of $\fol$.
\end{defi}

In Section~\ref{Sec:resolution}, another type of blow-up will be considered. 
In particular, in the course of Section~\ref{Sec:resolution}, blow-ups as in Definition~\ref{def:blowup}	will usually be referred to as {\it standard blow-ups}, while a different type of blow-ups will be called {\it weighted blow-up}. Let us provide an accurate definition for the latter.

\begin{defi}
	For $k \geq 2$, fix a $k$-tuple of strictly positive integers $\omega = (\omega_1, \ldots, \omega_k)$. The weighted projective space $\p^{k-1}_{\omega}$ associated with $\omega$ is a complex manifold of dimension $k-1$ defined as the quotient of $\C^k$ through the action of $\C^{\ast}$ defined as follows
	\[
	(\lambda, (y_1, \ldots, y_k)) \mapsto (\lambda^{\omega_1} y_1, \, \ldots, , \lambda^{\omega_k} y_k) \, .
	\]
	The vector $\omega = (\omega_1, \ldots, \omega_k)$ is called the weight vector.
\end{defi}

Consider then a manifold $M$ of dimension~$n$ and let $C$ be a submanifold of $M$. Assume that the submanifold $C$ is given, in certain local coordinates $(x_1, \ldots, x_n)$ for $M$, by $\{x_1 = \cdots = x_k = 0\}$. Fix then a vector $\omega = (\omega_1, \, \ldots, \, \omega_k)$.

\begin{defi}
	The blow-up of $M$ with center $C$ and weight $w$ is the submanifold of $\C^n \times \p^{k-1}_{\omega}$ defined as
	\[
	M_{C,\omega} = \{(x_1, \ldots, x_n, y_1, \ldots, y_k) \in \C^n \times \p^{k-1}_{\omega} \, : \, \,
	x_i^{\omega_j} \, y_j^{\omega_i} = \x_j^{\omega_i} \, y_i^{\omega_j}, \, \, \, 1 \leq i, j \leq k \} \, .
	\]
	The weighted blowing-up map is the restriction to $M_{C, \omega}$ of the natural projection
	\[
	pr: \C^n \times \p^{k-1}_{\omega} \to \C^n \, .
	\]
\end{defi}

To finish this section, let us introduce some terminology for $1$-dimensional foliations, that will be useful throughout the text. Denote by $\fol$ a $1$-dimensional foliation on a complex manifold $M$ and let $p \in M$ be a singular point of $\fol$. Next, let $X$ be a representative of $\fol$ on a neighborhood of $p$.

\begin{defi}
	The eigenvalues of $\fol$ at $p$ are the eigenvalues of $DX(p)$.
\end{defi}

We will say that the singular point $p$ is {\it elementary} if at least one of the eigenvalues of $\fol$ at $p$ is non-zero. Furthermore, the singular point $p$ is said to be a {\it saddle-node} if $\fol$ is elementary but has at least one eigenvalue equal to zero. The number of eigenvalues equal to zero is called the {\it rank} of the saddle-node.

We can also talk about the {\it rank of resonance} of $\fol$ at $p$ if its is applicable. To be more precise, let us recall what we mean by a {\it resonant singular} point.

\begin{defi}
	Let $\fol$ be a singular foliation and $p$ a singular point of $\fol$. Let $\lambda = (\lambda_1, \ldots , \lambda_n)$ be the vector of eigenvalues of $\fol$ at $p$. We say that $\fol$ is resonant at $p$ (or that the eigenvalues presents a resonant relation) if,
	for some $i$, there exists $I = (i_1, \ldots , i_n) \in \N_0^n$ with $\sum_{j=1}^n i_j \geq 2$ such that
	\[
	\lambda_i = (I, \lambda) = i_1 \, \lambda_1 + \ldots + i_n \, \lambda_n \, .
	\]
	If ${\rm dim} \, \{m \in \Z^n : (m, \lambda) = 0\} = k$, as vector space, then $\fol$ is called $k$-resonant.
\end{defi}

Finally, we will say that the $\fol$ (or, equivalently, the eigenvalues) is in the Siegel domain if the origin belongs to the convex hull of the eigenvalues in $\C$. Otherwise, we say that $\fol$ (or, equivalently, the eigenvalues) is in the Poincar\'e domain.


\section{Invariant analytic sets for (Lie algebras of) vector fields}\label{sec:invariant_sets}

\vspace{0.3cm}

The problem of existence of ``invariant manifolds'' has always been a central theme in the theory of dynamical systems.
Among others, these ``invariant manifolds'' usually provide reductions on the dimension of the corresponding phase-space.
For example, in the general theory of hyperbolic systems, the so-called stable manifolds are examples of invariant manifolds
and, in fact, their existence form a cornerstone of the hyperbolic theory.

In the local theory of vector fields a hyperbolic singular point of a vector field is an example of a hyperbolic set.
The existence of stable manifolds for such points is a consequence of the general theory and ensured by the well-known
Stable Manifold Theorem. Stable invariant manifolds, however, may fail to exist if the singular point is no longer
hyperbolic. For example, the integral curves associated to the (non-hyperbolic) vector field
\[
X = y \frac{\partial}{\partial x} - x \frac{\partial}{\partial y}
\]
are circles centered at the origin of $\R^2$ and the stable manifold is clearly empty in that case. Besides,
even in the case where the singular point is hyperbolic the stable manifold may not provide reduction on the
dimension of the corresponding phase-space. To find examples, it is sufficient to think of a planar vector
field whose hyperbolic singular point has two conjugated non-real eigenvalues, i.e. two non-real and non-pure
imaginary eigenvalues. In fact, in this case, the corresponding integral curves are spirals and the stable/unstable
manifold contains a neighborhood of the singular point.

The general problem of existence of ``invariant manifolds'' may also be considered in holomorphic dynamics. In this
case, however, there arise some important differences with the real counterpart. For example, in the holomorphic
setting we look for (proper) ``invariant manifolds'' that are analytic, which is a much stronger regularity condition.
We allow, however, the analytic manifolds to be singular in the sense of analytic sets, i.e. they are ``invariant
varieties'' as opposed to manifolds. In the sequel, the word ``manifold'' will be saved for smooth objects.

Briot and Bouquet were the first to consider the mentioned problem for holomorphic vector fields defined on a neighborhood
of the origin of $\C^2$. Basically they looked for the existence of separatrices (c.f. Definition~\ref{defi:separatrix}). In \cite{BB}, Briot and Bouquet claimed
the existence of separatrices for all holomorphic vector fields on $(\C^2,0)$. Their proof, however, contained a gap and their
classical work was completed only much later by Camacho and Sad. In fact, in their remarkable paper \cite{CS}, Camacho and Sad
prove the following.

\begin{teo}\cite{CS}\label{teo_CS}
	Let $\fol$ be a singular holomorphic foliation defined on a neighborhood of the origin of $\C^2$. Then there exists an analytic
	invariant curve passing through $(0,0)$ and invariant by $\fol$.
\end{teo}

Recall that in dimension~$2$ the singularities of every holomorphic foliation are necessarily isolated. Yet, the above Theorem applies
equally well to holomorphic vector fields. In fact, if $X$ is a holomorphic vector field on $(\C^2,0)$ with a curve of singular points,
then its components admit a non-invertible common factor, i.e. $X$ can be written as $X = fY$ for some holomorhic vector field with singular set of codimension at least~$2$ and where $f$ is a non-invertible holomorphic function. Up to eliminating this common factor $f$, we obtain the vector field $Y$ that is everywhere
parallel to $X$ and with isolated singular points. Theorem~\ref{teo_CS} can then be applied to $Y$ and this yields an invariant curve for $X$ as
well. Alternatively, even the curves of zero of $X$ might be thought of as an invariant curve for $X$, whether or not it is invariant for
the underlying foliation. In any case, the above argument applied to the vector field $Y$ shows that there always exist a curve invariant
by both $X$ and the underlying foliation.

It is surprising that separatrices for holomorphic vector fields on $(\C^2,0)$ always exist, despite the condition of analyticity
for the invariant curve. Note that the invariant curves for the holomorphic vector field $y \frac{\partial}{\partial x} - x
\frac{\partial}{\partial y}$ mentioned above are given by the two straight lines $y = \pm ix$, which are totally contained in
the non-real part of $\C^2$ (up to the singular point itself).

\begin{rem}
	{\rm In general, however, it is natural to allow the invariant curves
		to be singular at the singular point of the vector field, otherwise no general existence statement would hold. Indeed, as simple
		example, consider the holomorphic vector field $2y \partial /\partial x + x^3 \partial /\partial y$. Since this vector field admits
		$f(x,y) = x^3 - y^2$ as first integral, it immediately follows that the only separatrix of $X$ is the cusp of equation $\{x^3 - y^2
		= 0\}$, which is clearly not smooth at the origin. Fortunately, allowing separatrices to be singular is not a problem, as the classical
		theorem of resolution of singularities of Hironaka can be used to desingularize them, as well as general invariant analytic sets.}
\end{rem}

Unfortunately, the existence of separatrices is no longer a general phenomenon in dimension~$3$. To begin with, note that
when we move to dimension~$3$, it becomes necessary to distinguish between foliations of {\it dimension~$1$}\, and foliations
{\it of dimension~$2$ (or, equivalently, of codimension~$1$)}. Recalling Definition~\ref{defi:separatrix}, in the case of $1$-dimensional foliations, a separatrix is an (germ of)
analytic curve passing through the origin and invariant by the foliation while, in the context of codimension~$1$ foliations, a separatrix should be understood as a germ of a surface (i.e. a $2$-dimensional analytic set) passing through the
origin and invariant by the foliation. In fact, in dimension~$3$, the existence of separatrices is no longer a general
phenomenon, regardless of the dimension of the foliation, as it will made clear below.

Since for many years the works of Briot and Bouquet were thought to have established the general existence of separatrices in dimension~$2$, the question posed by Thom, on the existence of invariant varieties for codimension~$1$ foliations on $\C^3$, attracted much interest. A counterexample to his question was given by Jouanolou~\cite{Joa} on 1979 (slightly before Camacho and Sad completed Briot and Bouquet's work in dimension~$2$). Although the problem on the existence of separatrices in higher dimensions is {\it a priori} a local problem, as it will become clear in the course of the discussion below, the conterexample provided by Jouanolou is of {\it global nature}. The phenomena will be detailed below. In turn, examples of $1$-dimensional foliations without separatrices on $3$-manifolds were found by Gomez-Mont and
Luengo, \cite{GM-L}, some years later. Let us discuss both examples before providing an interesting ``generalization'' of Camacho and Sad Theorem under suitable conditions.

\bigbreak

\noindent{1. Gomez-Mont and Luengo counterexample for $1$-dimensional foliations}

\bigbreak

The example provided by Gomez-Mont and Luengo relies on a simple idea though its implementation requires significant computational
effort, which is carried with computer assistance. Yet, we may quickly describe the structure of their construction which relies on
two simple observations.

Consider then the foliation $\fol$ on $(\C^3,0)$ given by a holomorphic vector field satisfying the following conditions
\begin{itemize}
	\item[(1)] The origin $(0,0,0) \in \C^3$ is an isolated singularity of $X$
	
	\item[(2)] $J^1 X (0,0,0) = 0$ but $J^2 X (0,0,0) \ne 0$, where $J^k X (0,0,0)$ stands for the jet of order $k$ of $X$ at the origin
	($k=1,2$).
	
	\item[(3)] The quadratic part $X^2$ of $X$ at $(0,0,0)$ is a vector field whose singular set has codimension~$2$. In particular, $X^2$ is not
	a multiple of the Radial vector field $x \partial /\partial x + y \partial /\partial y + z \partial /\partial z$.
\end{itemize}
Assume that $\fol$ has a separatrix $C$ and consider the blow-up $\tilf$ of $\fol$ centered at the origin. Denote by $\pi$ the blow-up
map so that $\tilf = \pi^{\ast} \fol$ and let $\pi^{-1}(0)$ denote the exceptional divisor, which is isomorphic to $\C\p(2)$. Since, from item (3), $X^2$ is not
a multiple of the Radial vector field, there follows that $\pi^{-1}(0)$ is invariant by $\tilf$. Hence the restriction of $\tilf$ to
$\pi^{-1}(0)$ can be seen as a foliation of degree~$2$ on $\C\p(2)$.

Being $\pi^{-1}(0)$ invariant by $\tilf$, there follows that the transform $\pi^{-1}(C)$ of the separatrix $C$ can only intersect
$\pi^{-1}(0)$ at singular points of $\tilf|_{\pi^{-1}(0)}$. In other words, $\pi^{-1}(C)$ must be a separatrix (not contained in
$\pi^{-1}(0)$) for one of the singular points of $\tilf$.

Now, the second ingredient is as follows: being $\tilf|_{\pi^{-1}(0)}$ a foliation of degree~$2$ on $\C\p(2)$, $\tilf|_{\pi^{-1}(0)}$ has at most (and generically)
$7$ singular points. Since it is hard to control the position of $7$ points in $\C\p(2)$, the authors of~\cite{GM-L} proceed as follows.
\begin{itemize}
	\item[(A)] Let the singular points ``collide'' so as to have only~$3$ of them (position is then easily controlled)
	
	\item[(B)] Each of the~$3$ singular points will have an eigenvalue equal to zero associated to the direction transverse to $\pi^{-1}(0)$.
	The~$3$ singular points are therefore saddle-node singularities (in $3$-dimensional space).
	
	\item[(C)] Furthermore, arrange for the saddle-node singularities to have two equal eigenvalues tangent to $\pi^{-1}(0)$ that are, in addition, non-zero. The singular points in question are then (codimension~$1$) resonant saddle-nodes with weak direction transverse to $\pi^{-1}(0)$.
	
	\item[(D)] The three saddle-nodes are such that all separatrices are contained in $\pi^{-1}(0)$. In fact, it is well known that it is easy to produce saddle-node singular points with no separatrix not contained in the invariant $2$-plane associated with the non-zero eigenvalues.
\end{itemize}

The remainder of the proof of \cite{GM-L} consists of showing that it is, indeed, possible to prescribe a quadratic $X^2$ and a cubic $X^3$,
homogeneous components for the vector field $X$, so that all of the preceding conditions are satisfied.

Note that conditions (A), (B) and (C) depend only on the quadratic part $X^2$. The role played by the appropriated chosen cubic part $X^3$
can be summarized as follows.
\begin{itemize}
	\item it ensures each of the singular points of $\tilf$ are isolated singular points coinciding with the corresponding singular points of
	$\tilf|_{\pi^{-1}(0)}$. Here we note that the homogeneous foliation associated with $X^2$ has zeros all along the fibers of $\widetilde{\C}^3
	\to \widetilde{\pi}^{-1}(0)$ passing through the singular points of $\tilf|_{\pi^{-1}(0)}$. Thus some higher order perturbation
	to $X^2$ is already needed to have isolated singular points.

	\item having made sure the singular points are isolated, the cubic part $X^3$ of $X$ also takes care of condition (D)
\end{itemize}

As mentioned, the verification that all these conditions are compatible is conducted in \cite{GM-L} with the assistance of formal computations
programs.

\bigbreak

\noindent{2. Jouanoulu counterexample for codimension~$1$ foliations}

\bigbreak

The example of a codimension~$1$ foliation admitting no separatrix provided by Jouanolou is the foliation $\mathcal{D}_n$ defined as the kernel of the (integrable) $1$-form
\[
\Omega_n = (yx^n - z^{n+1}) \, dx + (zy^n - x^{n+1}) \, dy + (xz^n - y^{n+1}) \, dz \, ,
\]
with $n \in \N$. It can easily be checked that the kernel of $\Omega_n$ always contains the radial direction so that it naturally
induces a line field, and hence a ($1$-dimensional) foliation $\fol_n$, on $\C \p(2)$. The $1$-dimensional foliation $\fol_n$ can
also be viewed as follows: take the blow-up of $\C^3$ centered at the origin and let $\widetilde{\mathcal{D}}_n$ stands for
the transformed foliation of $\mathcal{D}_n$. Since the Radial vector field
\[
R = x \frac{\partial}{\partial x} + y \frac{\partial}{\partial y} + z \frac{\partial}{\partial z}
\]
is tangent to $\mcf_n$, the leaves of $\widetilde{\mcf}_n$ are (generically) transverse to the exceptional divisor which, in turn,
is isomorphic to $\C \p (2)$. The intersection of the leaves of $\tf_n$ with the exceptional divisor are then the leaves of the
$1$-dimensional foliation $\fol_n$ mentioned above.

The main result of Jouanolou states that $\fol_n$ leaves no algebraic curve invariant. This implies that $\mcf_n$ cannot admit a
separatrix. In fact, if $\mcf_n$ admits a separatrix $S$, i.e. an analytic surface $S$ that is invariant by $\mcf_n$, then the
intersection of the strict transform of $S$ by the mentioned blow-up at the origin with the exceptional divisor would be an invariant
algebraic curve for $\fol_n$, contradicting the result of Jouanolou.

\bigbreak

\noindent{3. Many other counterexample for codimension~$1$ foliations - a construction}

\bigbreak

There are many other examples of codimension~$1$ foliations on $\C^3$ without separatrix. Let me briefly explain how
numerous similar examples can be constructed. Consider a homogeneous polynomial vector field $Z$ defined on $\C^3$
and having an isolated singularity at $(0,0,0) \in \C^3$. Unless $Z$ is a multiple of the Radial vector field $R$,
it induces a $1$-dimensional holomorphic foliation $\fol$ on $\C \p (2)$. Conversely every $1$-dimensional foliation
on $\C \p(2)$ is induced by a homogeneous vector field on $\C^3$. Next, we consider the $2$-dimensional distribution
of planes on $\C^3$ that is spanned by $Z$ and $R$. The Euler relation (i.e. the equality $[R,Z] = (d-1)Z$, where $d$
is the degree of $Z$) shows that $Z,\, R$ generates a Lie algebra isomorphic to the Lie algebra of the affine group.
The corresponding distribution is therefore integrable and hence yields a codimension~$1$ foliation that will be
denoted by $\mcf$. Clearly the punctual blow-up $\tf$ of $\mcf$ does not leave the exceptional divisor $E \simeq \C \p(2)$
invariant (since the Radial vector field is tangent to $\mcf$) and thus $\tf$ induces a $1$-dimensional foliation on $E
\simeq \C \p(2)$. This $1$-dimensional foliation naturally corresponds to the intersection of the leaves of $\tf$ with
$E$. However, by construction, it also coincides with the leaves of $\fol$, the foliation induced by the homogeneous vector
field $Z$ on $\C\p(2)$. Note that the counterexample provided by Jouanolou fits this pattern. In fact, the Jouanolou foliation is tangent to both the Radial vector field $R$ and the homogeneous vector field $Z$ given by
\[
Z = y^n \frac{\partial}{\partial x} + z^n \frac{\partial}{\partial y} + x^n \frac{\partial}{\partial z} \, ,
\]
with $n \in \N$, $n \geq 2$.

As far as the existence of separatrices for $\mcf$ is concerned, the upshot of the preceding construction is as follows:
if $\mcf$ possesses a separatrix, the tangent cone of this separatrix yields an algebraic curve in $E \simeq \C \p(2)$
which must be invariant under $\fol$. Nonetheless, today it is known that, in a very strong sense, most choices of $Z$ lead to
a foliation $\fol$ that does not leave any proper analytic set invariant (cf. for example \cite{soares}, \cite{loray}).
As a result the codimension~$1$ foliation obtained by means of $Z, \, R$, for a generic choice of $Z$, does not have
separatrices. We also note that, for these examples, no separatrix can be produced by adding ``higher order terms'' to
$\mcf$. Jouanolou example fits this pattern.

\vspace{1cm}

This well known phenomena have led the experts (such as F. Cano, D. Cerveau and L. Stolovitch among others)
to wonder that the ``correct'' generalization of Camacho-Sad theorem would involve codimension~$1$ foliations
spanned by a pair of commuting vector fields (not everywhere parallel). The theorem below confirms their intuition
and affirmatively answers it. The proof can be found in~\cite{RR_separatrix}.

\begin{teo}\cite{RR_separatrix}\label{teo_RR_separatrix}
	Consider holomorphic vector fields $X, \, Y$ defined on a neighborhood of the origin of $\C^3$. Suppose that they
	commute and are linearly independent at generic points (so that they span a codimension~$1$ foliation denoted by
	$\mcf$). Then $\mcf$ possesses a separatrix.
\end{teo}

The existence of separatrices for codimension~$1$ foliations in general was also the object of some remarkable papers
such as \cite{CC} where it is proved the following

\begin{teo}\cite{CC}
	Let $\mathcal{D}$ be a non-dicritical codimension~$1$ foliation on a neighborhood of the origin of $\C^3$. Then $\mathcal{D}$ possesses a separatrix.
\end{teo}

However, as it follows from the preceding discussion, the set of foliations that fail to be {\it non-dicritical}\, is not negligible. In \cite{RR_separatrix} we establish the existence of separatrices for dicritical foliations that are spanned by two commuting vector fields.

The main difficulty in establishing the existence of a separatrix for a dicritical codimension~$1$ foliation lies in
controlling the dynamics of the $1$-dimensional foliations induced on the non-invariant, i.e. dicritical, components of
the exceptional divisor obtained after a suitable sequence of blow-ups. The key to prove the above theorem in our case was
to observe that these foliations always possess certain invariant algebraic curves provided that they are spanned by
commuting vector fields. It is the existence of these algebraic curves that leads to the existence of separatrices. As it
was to be expected, in the proof of our main result, it was needed to discuss the effect of the blow-up procedures of \cite{CC}
and \cite{cano} on the initial vector fields $X, \, Y$ and the fundamental desingularization results of these papers for
codimension~$1$ foliations played a role in our argument.

Let me be more precise concerning the idea of the proof of Theorem~\ref{teo_RR_separatrix}. Essentially, we have to show
that the phenomenon described above (i.e. in the case of a foliation generated by the Radial vector field and a homogeneous
holomorphic/meromorphic vector field) cannot take place in our context, unless the $1$-dimensional foliation induced on
$\C \p (2)$ admits certain invariant curves. To do that we shall consider the intersection of our codimension~$1$ foliation
$\mcf$ spanned by the commuting vector fields $X, Y$ with a given component of the exceptional divisor. Unless this component
is invariant by the codimension~$1$ foliation, this intersection defines a foliation of dimension~$1$ on it. Except for some
rather special situations that are already ``linear'' in a suitable sense, we are going to show that all the leaves of the
latter foliation are properly embedded (in particular they are compact provided that the mentioned component of the exceptional
divisor is so). This statement is, indeed, equivalent to saying that the corresponding foliation admits a non-constant meromorphic
first integral as it follows from \cite{Joa2}. In general we shall directly work with the existence of a non-constant meromorphic
first integral for foliations as above.

In view of the result proved by Cano and Cerveau, we have assumed the codimension~$1$ foliation $\mcf$ to be dicritical.
We assumed first that a non-irreducible component appears immediately after a single one-point blow-up along the assumption
that the first jet of $X$ and $Y$ at the point where we have centered the blow-up is zero. Since we are assuming the arising
exceptional divisor not to be invariant by the transformed foliation, we have that there exists a holomorphic vector field
$Z$ tangent to $\mcf$ and such that its first non-zero homogeneous component $Z^H$ of $Z$ is multiple of the Radial vector
field. There must then exist holomorphic functions $f, g$ and $h$ such that
\[
fX + gY = hZ \, .
\]
By exploiting the commutativity of the vector fields $X, Y$ and the assumptions of ``non-linearity'' of $X, Y$ we are able
to prove that the transformed foliation of $\mcf$ induces a $1$-dimensional foliation on the exceptional divisor with a
holomorphic/meromorphic first integral. In fact, it is the first non-trivial homogeneous components of $X, Y$ (denoted by
$X^H, Y^H$, respectively) that will play a role. Since $X,Y$ commute, so does $X^H, Y^H$. We can prove by using the above
relation that none of $X^H, Y^H$ is a multiple of the Radial vector field. There follows that both $X^H, \, Y^H$ induce a $1$-dimensional
foliation on $\C \p (2)$ by means of a punctual blow-up at the origin. Note, however, that the foliations induced by $X^H$ and by $Y^H$ must
coincide since $\C \p (2)$ is not invariant by $\tf$, the transform of $\mcf$. Furthermore, they must coincide with the
restriction of $\tf$ to $\C \p (2)$. It remains to prove that $X^H$ possesses a (non-constant) first integral and details
are provided in the paper.

Next we are led to consider the special situations of ``linear'' foliations that may not possess any non-constant first
integral. Fortunately, in these cases the existence of a separatrix can be established by more direct methods. An example of
a ``linear case'' would consist of a pair of vector fields $X, Y$ with $X$ linear and $Y$ equal to the Radial vector field. These
two vector fields commute and span a codimension~$1$ foliation whose (one-point) blow-up at the origin does not leave the corresponding
exceptional divisor invariant. Furthermore the foliation induced on the corresponding exceptional divisor by the mentioned blown-up
foliation coincides with the foliation induced on $\C \p (2)$ by $X$. In particular $X$ can be chosen so that the ``generic'' leaf
is not compact. However, in this situation the foliation induced by $X$ on $\C \p (2)$ still has a compact leaf which ``immediately''
leads to the existence of the desired separatrix.

Recall that the singular set of foliations on $3$-manifolds may have irreducible components of dimension~$1$. It was then necessary to derive an analogue of the above mentioned results for the case of blow-ups centered at a (smooth and irreducible) curve. Also in this case we have considered separately the ``linear'' and ``non-linear'' case and the following was proved: in the ``non-linear'' case the foliation induced on the exceptional divisor admits a non-constant first integral while in the ``linear'' case the existence of at least a contact leaf for the foliation on the exceptional divisor is proved. The words ``linear'' and ``non-linear'' are between quotes since we had to adapt the notion of ``linearity'' in
this case. To explain the need for this adaption (that will be made explicit below), let me describe an example that
was pointed out to us by D. Cerveau and that illustrates the problem about the existence of first integrals as above as also
some intermediary results which are crucial for establishing the existence of these first integrals in the case that a blow-up
along a (smooth) curve is considered. This example goes as follows. Consider the pair of vector fields $X, Y$ given by
$$
X = zy \frac{\partial}{\partial y} + z^2 \frac{\partial}{\partial z} \; \qquad {\rm and} \; \qquad Y = x^2 \frac{\partial}{\partial x} + axy \frac{\partial}{\partial y}
$$
These two vector fields commute and span a codimension~$2$ foliation
denoted by $\mcf$. They also leave the axis $\{ y=z=0\}$ invariant.
Consider the blow-up of $\mcf, X, Y$ centered over $\{ y=z=0\}$. The
transform $\tf$ of $\mcf$ does not leave the exceptional
divisor invariant. Furthermore the leaves of the foliation
induced on the non-compact exceptional divisor by intersecting it with the
leaves of $\tf$ are themselves non-compact. The explanation for this phenomenon
is that the blow-up of $Y$ is regular at generic points of the
exceptional divisor. Indeed, $Y$ is already regular at generic
points of the axis $\{ y=z=0\}$ - although quadratic (in the usual sense) at the origin. Hence this case must be considered
as ``linear'' (indeed even ``regular''). It then follows that the appropriate
notion of order of a vector field relative to a curve is such that the resulting
order for $X$ as above is {\it zero}. In fact, if we proceed to blow-up
the curve locally given by $\{y = z = 0\}$, we must consider the variable $x$ as a constant. In this sense, the vector field is considered quadratic while the vector field $Y$ is considered regular. The reader is referred to~\cite{RR_separatrix} for a
precise definition of what is meant by ``linear'' or ``non-linear'' in this case.

\bigbreak

Theorem~\ref{teo_RR_separatrix} established the existence of separatrices for the codimension~$1$ foliations spanned
by two commuting vector fields. We can ask what happens to the vector fields $X,Y$ themselves or to the foliations induced
by them. In this context, Theorem~\ref{teo_RR_separatrix} was complemented by another result in~\cite{RR_secondjet}. More
precisely, the following has been proved.

\begin{teo}\cite{RR_secondjet}\label{teo_RR_secondjet}
	Let $X$ and $Y$ be two holomorphic vector fields defined on a neighborhood $U$ of $(0,0,0) \in \C^3$ which are not linearly
	dependent on all of $U$. Suppose that $X$ and $Y$ vanish at the origin and that one of the following conditions holds:
	\begin{itemize}
		\item $[X,Y]=0$;
		\item $[X,Y]=c \, Y$, for a certain $c \in \C^{\ast}$.
	\end{itemize}
	Then there exists a germ of analytic curve $\mathcal{C} \subset \C^3$ passing through the origin and simultaneously
	invariant under $X$ and $Y$.
\end{teo}

This result deserves some comments. First of all, we should note that the Theorem~\ref{teo_RR_secondjet} applies not only to the commutative Lie algebras but also to the Lie
algebra of the affine group. Recall, however, that the analogue of Theorem~\ref{teo_RR_separatrix} in the case of affine
actions is known to be false since the classical work of Jouanolou. Whereas Theorem~\ref{teo_RR_secondjet} holds interest
in its own right as a theorem claiming the existence of invariant manifolds (or, more precisely, curves in this case), our paper~\cite{RR}
also contains a non-trivial application of this result (cf. Theorem~\ref{teo_RR_nonvanishing_secondjet} in Section~\ref{sec:localaspects}).

Secondly, Theorem~\ref{teo_RR_secondjet} states that $X$ and $Y$ possess a common invariant curve without mentioning if
the curve in question is invariant for the associated foliations (recall, for example, that $\{x=0\}$ is invariant by $x
\partial /\partial x$ but not by the corresponding foliation). It is however easy to check that in the particular case
that we consider $X$ as being a homogeneous vector field and $Y$ as a multiple of the Radial vector field, the existence
of a common separatrix for $\fol_X$ and $\fol_Y$ can easily be deduced. In fact, the leaves of $\fol_Y$ are simply the
radial lines. Concerning $\fol_X$, since it is not a multiple of the Radial vector field, it induces a $1$-dimensional
foliation on $\C \p (2)$ by means of the one-point blow-up of $\C^3$ at the origin. The foliation in question possesses
isolated singular points and it can easily be checked that the radial line naturally associated with any of these singular
points is invariant by $\fol_X$ as well. We believe that the existence of a common separatrix for $\fol_X$ and $\fol_Y$
in the general case can also be established.

To finish this section we are just going to give an idea of the proof of Theorem~\ref{teo_RR_secondjet}.
Let then $\mcf$ stands for the codimension~$1$ foliation spanned by $X$ and $Y$. We have that ${\rm codim} \, ({\rm Sing} \,
(\mcf)) \geq 2$. In other words, ${\rm Sing} \, (\mcf)$ is of one of the following types: 
\begin{itemize}
	\item the union of a finite number of irreducible curves,
	\item a single point (the origin) or
	\item empty (which means that $\mcf$ is regular).
\end{itemize}

Since ${\rm Sing} \, (\mcf)$ is naturally invariant by $X$ and $Y$, the result immediately holds if ${\rm dim} \, ({\rm Sing} \, (\mcf))=1$. Hence we can assume without loss of generality that ${\rm Sing} \, (\mcf)$ has codimension at least~$3$.

Since the singular set of $\mcf$ has codimension at least~$3$, $\mcf$ possesses a holomorphic first integral $f$ thanks to Malgrange
Theorem (cf. \cite{malgrange}). Let then $S = f^{-1}(0)$. We have that $S$ is an invariant surface for $\mcf$ and, consequently, for $X$ and $Y$. Furthermore,
$S$ can be assumed to be irreducible. In fact, if it was not, then the intersection of any two irreducible components of $S$ is an
invariant curve for both $X, \, Y$ and the conclusion holds. The surface $S$ can be assumed either regular or having an isolated singularity
at the origin (in fact, if the singular set of $S$ contains a singular curve, then the singular curve is a common separatrix for $X, \, Y$).
Next, the following can be noted
\begin{itemize}
	\item In the case where $S$ is smooth, both foliations ($\fol_X$ and $\fol_Y$) possess a separatrix through the origin by Camacho-Sad Theorem. We have
	however to check that the imposed conditions implies that at least one of their separatrices coincide.
	
	\item In the case of singular surfaces, there are examples of holomorphic vector fields without separatrix (cf.~\cite{C}).
	This phenomenon needs to be ruled out in the case in question.
\end{itemize}

Consider then the restrictions of $X$ and $Y$ to the invariant surface $S$ along with the corresponding tangency locus (which is non-empty since the origin
is a common singular point of $X, \, Y$). Since the tangency locus ${\rm Tang} \, (X|_S, Y|_S)$ is invariant by both $X$ and $Y$, the
result immediately holds in the case where its dimension equals 1. So, we shall consider separately the case where ${\rm Tang} \, (X|_S,
Y|_S) = \{0\}$ and ${\rm Tang} \, (X|_S, Y|_S) = S$

Assuming that ${\rm Tang} \, (X|_S, Y|_S) = \{0\}$, we get that $S$ is a surface with singular set of codimension at least 2 and
equipped with two linearly independent vector fields. This implies that the tangent sheaf to $S$ is locally trivial which, in turn, implies
that $S$ is smooth. However, being $S$ smooth, we have that $S$ is locally equivalent to $C^2$ and the tangency locus of two vector
fields on there cannot be reduced to a single point. The contradiction excludes this case.

We should then assume that ${\rm Tang} \, (X|_S, Y|_S) = S$, i.e. $X$ and $Y$ coincide up to a multiplicative function on $S$. The
existence of the desired common separatrix is then ensured in the case where $S$ is smooth. It remains to consider the case where $S$
is singular at the origin. The argument in this case relies on proving that the ($1$-dimensional) foliation induced on $S$ by either
$X$ or $Y$ possesses a non-constant holomorphic first integral. The level curve of this first integral containing the origin then
yields the desired separatrix. Details can be found in the paper in question.


\section{Vector fields with univalued solutions - Global dynamics}\label{sec:sc_global_dynamics}

\begin{figure}[hbtp]
	\centering
	\includegraphics[scale=0.7]{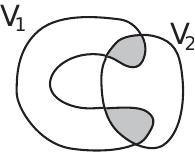}
	\caption{Domains with non-connected intersection}
	\label{fig:monodromy_effect}
\end{figure}

Recall that a solution $\varphi$ of a real ordinary differential equation (ODE) always possesses a maximal domain of definition contained in $\R$. In other
words, fixed an initial condition, the solution $\varphi$ is defined on a maximal interval $I_0 \subseteq \R$ in the sense that if $\{t_i\}
\subseteq I_0$ is a sequence converging to an endpoint of $I_0 (\ne \pm \infty)$, then the sequence $\{\varphi(t_i)\}$ leaves every compact
set in $M$ as $i \to \infty$. Unlike the real case, the solution of a complex ODE (i.e. those where the time is a complex variable) does not
admit in general a maximal domain of definition contained in $\C$. In fact, those admitting a maximal domain of definition are somehow ``rare''
among all complex ODE's. The absence of maximal domains of definition for solutions of complex ODE's is closely related to standard ``monodromy''
phenomena arising when extending holomorphic functions along paths. This phenomenon is illustrated by Figure~\ref{fig:monodromy_effect}, since the intersection of the domains $V_1$ and $V_2$ is not connected, solutions defined on $V_1$ and $V_2$ cannot, in general, be ``adjusted'' to coincide in both connected components. In view of this, complex ODE whose solutions admit maximal domains of
definition are, roughly speaking, those whose solutions are univalued.

The understanding of the mentioned equations, which is a far classical and important problem with many interesting applications,
is one of the topics in my research. It was probably Painlev\'e who first paid attention to these problems while studying equations
not admitting movable singular points. Painlev\'e's motivations were mostly concerned with the theory of special functions
and they remain an active area of research nowadays. These and other problems connected to special function theory also
constitute a motivation for my own past and future work.

On the other hand, in algebraic/complex geometry, there is also a fundamental problem of describing the ``pairs'' consisting
of a holomorphic vector field defined on compact manifold. Since we are dealing with compact manifolds, vector fields as above
are necessarily {\it complete}\, in the sense that their solutions are defined on all of $\C$. In some more specific cases (e.g.
affine geometry, groups of birational automorphisms etc), one also pays attention to the problem of classifying complete holomorphic
vector fields (on open manifolds). In both cases, a wealth of information is encoded in the nature of the singular points of the
corresponding vector fields.

A surprising connection between the study of the above mentioned singular points and the existence of maximal
domains of definitions for solutions of complex ODEs was realized by Rebelo, who introduced the notion of {\it semicomplete singularity}\,
in \cite{Rebelo96}. The idea is that the solutions of these (local) vector fields must be univalued since they have a ``realization''
as a complete vector field on some complex manifold. The condition of being ``semicomplete'' for a (germ of) vector field turned out
to be very non-trivial and, indeed, to capture almost all of the ``intrinsic nature'' of germs of vector fields that actually admit
global realizations as complete vector fields. The ``classification'' of semicomplete singularities of vector fields has then become
an important problem which, a {\it posteriori}, has also shown a number of interesting connections with integrable systems and certain
remarkable kleinian groups, as will be seen later.

To begin with, let us introduce some definitions and results that will be useful in the sequel. So, let $X$ be a holomorphic vector
field defined on a possibly open complex manifold $M$ and let $U$ be an open subset of $M$.

\begin{defi}
	The holomorphic vector field~$X$ is said to be semicomplete on $U$ if for every~$p\in U$ there exists a connected domain~$V_p
	\subset \C$ with~$0\in V_p$ and a map~$\phi_p : V_p \to U$ satisfying the following conditions:
	\begin{itemize}
		\item $\phi_p(0)=p$
		\item $\phi'_p(T) = X(\phi_p(T))$, for every $T \in V_p$.
		\item For every sequence~$\{T_i\}\subset V_p$ such that~$\lim_{i\rightarrow\infty} T_i = \hat{T} \in \partial V_p$ the sequence
		$\{\phi_p(T_i)\}$ escapes from every compact subset of~$U$.
	\end{itemize}
\end{defi}

In this way $\phi : V_p \subset \C \rightarrow U$ is a maximal solution of $X$ in a sense similar to the notion of ``maximal solutions''
commonly used for real differential equations. A semicomplete vector field on $U$ gives rise to a semi-global flow $\Phi$ on $U$.

It should be mentioned that if $X$ is semicomplete on $U$ and $V \subset U$, the restriction of $X$ to $V$ is semicomplete as
well. Thus the notion of ``semicomplete singularity'' (or germ of semicomplete vector field) is well defined. It immediately
follows that if $X$ is globally defined on a compact manifold $M$ then $X$ is semicomplete at every singular point. In this
sense, semicomplete vector fields can be viewed as the ``local version'' of complete ones. In fact, a singularity that is not
semicomplete cannot be realized by a complete vector field. In particular, it cannot be realized by a globally defined
holomorphic vector field on a compact manifold. Yet, the same definition applies also to more global context since the
set $U$ need not be ``small''. This is especially meaningful in the context of rational/polynomial vector fields that may
be semicomplete away from their pole divisors. In these situation we shall use the terminology {\it uniformizable vector
	field} so as to save the phrase ``semicomplete vector field'' for situations where we shall be working on a neighborhood
of a singular point.

A useful criterion to detect semicomplete vector fields was deduced by Rebelo in \cite{Rebelo96} and can be stated as follows. First consider a holomorphic vector field $X$ on $U$
and note that the local orbits of $X$ define a singular foliation $\fol$ on $U$. A regular leaf $L$ of $\fol$ is naturally a Riemann
surface equipped with an Abelian $1$-form $dT_L$ which is called the {\it time-form}\, induced on $L$ by $X$. Indeed, at a point $p
\in L$ where $X(p) \neq 0$, $dT_L$ is defined by setting $dT_L (p).X(p) =1$. Now, we have the following

\begin{prop}\cite{Rebelo96}
	Let $X$ be a holomorphic vector field defined on a complex manifold $M$. Assume that $X$ is semicomplete on $M$. Then
	\[
	\int_c dT_L \ne 0
	\]
	for every open (embedded) path $c: [0,1] \rightarrow L$.
\end{prop}

This proposition allows us to easily check that the vector field $X = x^3 \partial /\partial x$, for example, is not semicomplete on any neighborhood
of the origin of $\C$. In fact, for every given neighborhood $U$ of the origin, there exists $\varepsilon > 0$ such that the
ball $B_{\varepsilon}$ with center at the origin and radius $\varepsilon$ is contained in $U$. Consider then the open path contained in
$U$ and given by $c(t) = \frac{\varepsilon}{2} e^{\pi i t}$, $0 \leq t \leq 1$. Since the time-form on the (unique) regular leaf associated
to $X$ is given by
\[
dT = \frac{dx}{x^3} \, ,
\]
it becomes clear that the integral of the time form along the open path in question is zero and the claim immediately follows. We can
also check that the vector field in question is not semicomplete on any neighborhood of the origin by noticing that the solution of the
differential equation associated to it
\[
x(T) = \frac{x_0}{\sqrt{1 - 2x_0^2 T}}
\]
is multivalued, where $x_0 = x(0)$.

In turn, it can easily be checked that the vector field $X = x^2 \partial /\partial x$ is semicomplete on a neighbohood of the origin. In fact, it is semicomplete on the whole of $\C$. To check this, it suffices to notice that the solution of the differential equation associated with $X$ satisfying $x(0)=x_0$ takes on the form
\[
x(t) = \frac{x_0}{1-x_0 \, T} \, .
\]
The solution is defined on $V = \C \setminus \{1/x_0\}$ and it ``escapes'' to infinity whenever $T$ goes to $1/x_0$, which belongs to the boundary of $V$. However, an alternative proof that $X$ is semicomplete on $\C$ can be obtained by observing that $X$ admits a holomorphic extention to $\C\p(1)$ and by applying the results above.

In terms of classification of semicomplete vector fields on $(\C,0)$ we have the following

\begin{prop}\cite{Rebelo96}
	If $X$ is a semicomplete (holomorphic) vector field on a neighborhood of $0 \in \C$ and if its first jet at the origin vanishes identically, then $X$ is analytically conjugated to the vector field $x^2 \partial/\partial x$.
\end{prop}

Recall that

\begin{defi}
	Let $X, \, Y$ be two vector fields defined on a neighborhood of the origin of $\C^n$. The vector fields are said to be analytically conjugated if there exists a holomorphic diffeomorphism $H: (\C^n,0) \to (\C^n,0)$ such that $DH \, . \, Y = X \circ H$. In the case where $DH \, . Y = f \, . \, (X \circ H)$, for some holomorphic function $f$, then $X$ and $Y$ are said analytically equivalent.
\end{defi}

\begin{rem}
	{\rm If two vector fields are analytically conjugate and one of them is semicomplete, then so is the other. The same does not hold for analytically equivalent vector fields, as the previous example makes it clear. It should also be mentioned that the semicomplete character of a vector field is preserved by biratioal transformations, in particular by blow-ups.}
\end{rem}

In dimension~2, semicomplete singularities of vector fields (whether or not isolated) were fully classified by Ghys and Rebelo
and this classification was, in particular, strongly used in the description of pairs of compact complex surfaces equipped with
a globally defined holomorphic vector field obtained in~\cite{dloussky}. In terms of semicomplete vector fields with an isolated singularity the following can be said.

\begin{teo}\cite{RG}\label{teo:sc_quadratic}
	Let $X$ be a holomorhic vector field defined on a neighborhood of the origin of $\C^2$ and such that the origin is an isolated singularity for $X$. Assume in addition that the first jet of $X$ at the origin vanishes identically. If $X$ is semicomplete in a neighborhood $U$ of the origin, then $X$ is analytically conjugate to one of the vector fields:
	\begin{itemize}
		\item[1.] $f \, [x \partial/\partial x - y(nx - (n + l)y) \partial/\partial y]$, where $n$ is a strictly positive integer;
		
		\item[2.] $f \, [x(x - 2y) \partial/\partial x + y(y - 2x) \partial/\partial y]$;
		
		\item[3.] $f \, [x(x - 3y) \partial/\partial x + y(y - 3x) \partial/\partial y]$;
		
		\item[4.] $f \, [x(2x - 5y) \partial/\partial x + y(y - ^x) \partial/\partial y]$,
	\end{itemize}
	where $f$ is a holomorphic function on $U$ with $f(0,0) \ne 0$.
\end{teo}

They have proved that a semicomplete vector field as in the above theorem needs to be analytically equivalent to a quadratic homogeneous vector field, although not necessarily analytically conjugated. This issue, however, is compensated by the use of the multiplicative function``$f$'' As it will be recalled later, the vector fields in question are integrable. In fact, the first one admits a meromorphic first integral while the other three possess a holomorphic first integral.

The extension of their results to higher dimensions
is, however, a very challenging and wide open problem, already in dimension 3. In fact, the study of semicomplete vector
fields in dimensions~$\geq 3$ was initiated by A. Guillot \cite{guillotFourier}, \cite{guillotIHES} who sought to classify
quadratic semicomplete vector fields on $(\C^3,0)$ (since the vector fields are homogeneous, they are semicomplete on a
neighborhood of the origin if and only if they are ``uniformizable'' on all of $\C^n$, c.f. Corollary~2.6 in \cite{RG}). The interest in homogeneous vector
fields comes, in part, from the fact that semicompleteness is closed for the topology of uniform convergency on compact
subsets and, consequently, if a given vector field is semicomplete then so is its first non-zero homogeneous components.
Another motivation for Guillot's work stemmed from the evidence that among semicomplete vector fields one often finds
especially interesting/remarkable examples of dynamical systems, an idea totally in line with Painlev\'e's point of view
concerning equations without movable singular points. It is fair to say that some of the main outcomes in Guillot's work
concern the description of certain examples exhibiting remarkable properties in a way or in another.

Some of my works are contributions to the study of uniformizable vector fields on higher dimensional manifolds
(cf. \cite{Reis06}, \cite{Reis08}, \cite{RR_applications}, \cite{RR_secondjet}). The results with global nature will be discussed
below, while the results with local nature will be discussed in the next section.

Let us focus on the class of uniformizable vector fields from a definitely global point of view: the vector fields in
question are polynomial on $\C^n$ and are supposed to be semicomplete on all of $\C^n$ (as above mentioned, given the
global nature of the discussion we shall use the terminology ``uniformizable'' instead of ``semicomplete''). As it will
be pointed out later, the methods to be discussed in this section apply also to uniformizable rational vector fields
(where uniformizable means semicomplete away from its pole divisor). Examples fitting in this context include
complete polynomial vector fields but also certain uniformizable vector fields with solutions defined on hyperbolic
domains (some of them defined on a disc) as it happens in the case of Halphen vector fields. The works of Ablowitz
and his collaborators on evolutions equations - many of them appearing in fluid dynamics - also provides numerous examples
of equations/vector fields to which our methods are applicable, see for example \cite{ablowitz1and2}
and references therein.

This said, in the paper \cite{RR_applications}, a method to investigate the domain of definition of solutions for polynomial
(or, more generally, rational) vector fields was introduced. The method is quite general in that it applies to arbitrarily high
dimensions. Yet, it provides new results already in dimension~$3$. The mentioned paper fundamentally consists of two parts,
the first one corresponding to a general setup along with the basic estimates/results whereas the second part provides some
applications of it. To greater or lesser extent, the applications given there arise from following the solution of a (complex)
polynomial/rational vector field over ``special real paths going off to infinity''. Recall that being the vector fields
considered polynomial vector fields in $\C^n$, they admit a meromorphic extension to $\C \p (n)$. In particular, they induce
a holomorphic foliation $\fol$ in $\C \p (n)$. Let us denote by $\Delta_{\infty}$ the hyperplane at infinity of $\C \p (n)$,
i.e. let $\Delta_{\infty} = \C \p (n) \setminus \C^n$. We have that $\Delta_{\infty}$ is invariant by $\fol$ unless the top
degree homogeneous component of $X$ is a multiple of the Radial vector field.

Before describing the method in question, let us present some of the applications of it. The first application obtained concerns a
confinement-type theorem for solutions of complete polynomial vector fields on $\C^n$.

\begin{teo}\cite{RR_applications}\label{teo_area}
	Suppose that $X$ is a complete polynomial vector field of degree at least~$2$ on $\C^n$. Fix an arbitrarily small neighborhood
	$V$ of $({\rm Sing}\, (\fol) \cap \Delta_{\infty}) \cup {\rm Sing}\, (X)$ in $\C \p(n)$ and suppose we are given a point $p \in
	\C^n$, $X (p) \neq 0$, and an angle $\theta \in (-\pi/2 ,\pi/2)$. Denote by $L_p$ the leaf of $\fol$ through $p$ and consider the
	parametrization of $L_p$ by $\C$ (possibly as a covering map) which is given by $\Phi_p (T) = \Phi (T, p)$. Then there exists a
	separating curve $c: (-\infty, \infty) \rightarrow \C$, $\Phi_p(c(0)) =p$, and an unbounded component $\mathcal{U}^+$ of $\C
	\setminus \{c(t)\}$ such that the following holds: the set $\mathcal{T}_V \subset \mathcal{U}^+ \subset \C$ defined by
	$$
	\mathcal{T}_V = \{ T \in \mathcal{U}^+ \subset \C \; : \; \Phi (T,p) \in V \} \,
	$$
	satisfies
	$$
	\lim_{r \rightarrow \infty} \frac{ {\rm Meas}\, (\mathcal{T}_V \cap B_r)}
	{{\rm Meas}\, (\mathcal{U}^+ \cap B_r)} =1 \, ,
	$$
	where ${\rm Meas}$ stands for the usual Lebesgue measure of $\C \, (\simeq \R^2)$ and $B_r$ the ball of radius $r$ centered at $0 \in \C$.
\end{teo}

\begin{figure}[hbtp]
	\centering
	\includegraphics[scale=0.75]{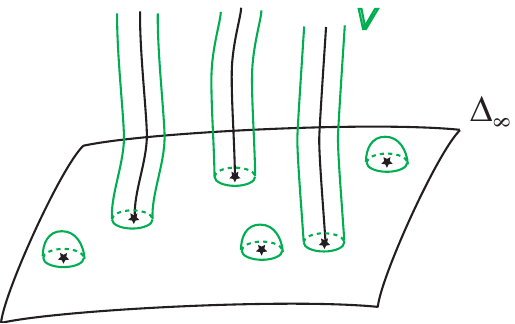}
	\caption{A neighborhood $V$ of the singular set}
	\label{fig:singularset}
\end{figure}

The ``separating curve'' is actually a geodesic for some singular flat structure on $\C$ having bounded coefficients with respect
to the standard flat structure. As a consequence of this fact, it follows that ${\rm Meas}\, (\mathcal{U}^+ \cap B_r)$ is actually
comparable to the measure of the large discs $B_r$. This can naturally be seen as a confinement phenomenon since the solutions of
a complete polynomial vector field spend a significant ``part of their existence in arbitrarily small regions of the phase space''
and hence ``are highly non-ergodic''. This phenomenon of ``strong non-ergodicity'' becomes more clear after the following corollary
of the preceding theorem.

\begin{corollary}\cite{RR_applications}
	Let us keep the notations of Theorem~\ref{teo_area} and let $\mathcal{T}_V(r)$ be the set of
	\[
	\mathcal{T}_V (r) = \{ T \in B_r \subset \C \; ; \; \Phi (T,p) \in V \} \, .
	\]
	Then, there exists $\delta > 0$ uniform (i.e. not depending on the neighborhood $V$) such that
	\[
	\liminf_{r \rightarrow \infty} \frac{ {\rm Meas}\, (\mathcal{T}_V (r))}
	{{\rm Meas}\, (B_r)} \geq \delta >0 \, .
	\]
\end{corollary}

Since the above inequality remains valid if we reduce~$V$, we see that the frequency with which the point~$p$ visits the neighborhood~$V$
is ``far'' from being proportional to the size of $V$.

The statement of Theorem~\ref{teo_area} and of the corresponding corollary indicate that the structure of the singularities of $\fol$
lying in $\Delta_{\infty}$ must hold significant information on the global dynamics of corresponding vector fields. This is a principle
similar to the ``heuristics'' involved in the classical Painlev\'e test: the local behavior at singular points have strong influence on
the global dynamics of the system. To put this principle to test, we have next assumed that the singularities of $\fol$ lying in $\Delta_{\infty}$
are ``simple'' (see below for the precised definition of ``simple'' singularity). The idea is that these singularities in particular can be understood in detail and
hence we should be able to derive strong consequences of the global behavior of the corresponding polynomial vector field.
Theorem~\ref{teo_simple_sing} then fully vindicates our principle. Before stating the theorem, let us just mention what we mean by ``simple'' singularity. For us, a singularity $q \in \Delta_{\infty}$ for $\fol$ is said {it simple} if it is of one of the following types:
\begin{itemize}
	\item[(1)] A non-degenerate singularity: this means that $\fol$ can locally be represented by a vector field having non-degenerate
	linear part at $q$ (i.e. the Jacobian matrix of $X$ at $q$ is invertible or, equivalently, it possesses~$n$ eigenvalues different from
	zero). Besides, since resonances may arise, we assume that $q$ is not of Poincar\'e-Dulac type, i.e. if all the eigenvalues of $\fol$
	at $q$ belong to $\R^{\ast}_+$, then $\fol$ must be locally linearizable about $q$.
	
	\item[(2)] Codimension~$1$ saddle-node: these are singularities of $\fol$ lying in $\Delta_{\infty}$ whose eigenvalue associated
	to the direction transverse to $\Delta_{\infty}$ is equal to zero whereas it has $n-1$ eigenvalues different from zero and corresponding
	to directions contained in $\Delta_{\infty}$. Again we require that the singularity for the $(n-1)$-dimensional foliation induced on the
	plane $\Delta_{\infty}$ should not be a singularity of Poincar\'e-Dulac type.
\end{itemize}

Thus we have the following.

\begin{teo}\cite{RR_applications}\label{teo_simple_sing}
	Let $X$ be a complete polynomial vector field on $\C^n$ whose singular set has codimension at least~$2$. Suppose that all singularities
	of $\fol$ lying in $\Delta_{\infty}$ are simple. Then all leaves of $\fol$ can be compactified into a rational curve, i.e. $\fol$ can be
	pictured as a ``non-linear rational pencil''.
\end{teo}

Before proceeding and mentioning another application of these techniques (also provided in \cite{RR_applications}), let me describe the
method introduced in the mentioned paper along with the main ideas for the proofs of the theorems stated above.

So, let $X$ be a polynomial semicomplete vector field on all of $\C^n$.
Our method relies on estimating the ``speed'' of the vector field $X$ near $\Delta_{\infty}$, the hyperplane at infinity
(which coincides with the divisor of poles since $X$ is polynomial). This is done in two steps. The first step consists of
eliminating the unbounded factor of $X$ over $\Delta_{\infty}$ so as to obtain a ``local regular vector field'' about every
regular point $p \in \Delta_{\infty}$ of $\fol_{\infty}$, where $\fol_{\infty}$ stands for the foliation induced by $X$ at
$\Delta_{\infty}$. Recall that being $X$ polynomial, it admits a meromorphic extension to the plane at infinity inducing,
in particular, a holomorphic foliation on that. However, it turns out that these locally defined vector fields, obtained by
eliminating the unbounded factor, depend to some extent on the choice of local coordinates so that they do not patch together
in a ``foliated'' global vector field. Nonetheless, two local representatives obtained through overlapping coordinates differ
only by a (non-zero) multiplicative constant. This means that this collection of local vector fields defines a global affine structure
on every leaf of $\fol_{\infty}$. The interest on the mentioned affine structure lies in the fact that it lends itself well
to provide estimates for the flow of $X$ as long as accurate estimates for the ``distance'' from the orbit in question to
$\Delta_{\infty}$ are available.

\begin{rem}
	Recall that an affine (resp. translation) structure on a Riemann surface $L$ is nothing but a collection of charts on $L$ such that the change of charts are restrictions of affine (resp. translation) maps of $\C$. In the case where $L$ is a leaf of $\fol$ that is not contained in $\Delta_0 \cup \Delta_{\infty}$, the restriction of $X$ to $L$ is neither identically zero nor identically infinity. More precisely, $L$ becomes globally equipped with a non-identically zero holomorphic vector field and the local solutions of the differential equation associated to it naturally induce a translation structure on $L$. The same does not occur with respect to a leaf $L_{\infty}$ of $\fol$ on $\Delta_{\infty}$. In turn, the translation structures on the leaves nearby $L_{\infty}$ that are not contained in $\Delta_{\infty}$ will induce an affine structure on $L_{\infty}$.
\end{rem}

The second ingredient of our construction is a quantitative measure of ``the rate of approximation'' of a leaf of $\fol$
to $\Delta_{\infty}$. Because $\Delta_{\infty} \subset \C P(n)$ and the Fubini-Study metric on $\C P(n)$ has positive
curvature, it is well known that complex submanifolds always ``bend themselves towards $\Delta_{\infty}$''. In our case,
this implies that the distance (relative to the Fubini-Study metric) of a leaf $L$ of $\fol$ to $\Delta_{\infty}$ can never
reach a local minimum unless this minimum is zero. Our mentioned second ingredient is reminiscent of this remark though,
in the mentioned paper, the euclidean metric on suitably chosen affine coordinates, as opposed to the globally defined
Fubini-Study metric, was chosen. The choice is however a relatively minor technical point due to the fact that the
euclidean metric is better adapted to work with the above mentioned affine structure. Besides,
by exploiting the fact that the submanifolds in questions are actual leaves of a foliation, a quantitative version
of the rate of approximation of a leaf to $\Delta_{\infty}$ is derived. The phenomenon goes essentially as follows.
At each regular point $p$ of a leaf $L$ of $\fol$ there is
the ``steepest descent direction of $L$ towards $\Delta_{\infty}$'', namely the negative of the gradient of the
distance function restricted to $L$. This yields a singular real one-dimensional oriented foliation $\calh$ on $L$.
Roughly speaking, an exponential rate of approximation for $L$ to $\Delta_{\infty}$ over the trajectories of $\calh$ can then
be obtained. Since $L$ is endowed with a conformal structure, it makes sense to define also foliations $\calh^{\theta}$
whose (oriented) trajectories makes an angle~$\theta$ with the oriented trajectories of
$\calh$ ($\theta \in [-\pi /2, \pi /2]$). For $\theta \in ]-\pi /2, \pi /2[$ an exponential rate of approximation for
$L$ to $\Delta_{\infty}$ over the trajectories of the associated real foliation can also be obtained (note that the
foliation $\calh^{\pi /2}$, which is orthogonal to $\calh$, is constituted by level curves for the
above mentioned distance function). Finally, the estimates on the exponential rate of approximation combines to
the ``uniform'' estimates related to the foliated affine structure to produce accurate estimates for the time taken
by the flow of $X$ over trajectories of $\calh$.

To better explain the method, assume for simplicity that $X$ is a (polynomial) homogeneous vector field of degree $d \geq 2$
on $\C^3$. Assume, in addition, that $X$ is not a multiple of the Radial vector field. Let us consider $\C \p (3)$, the
compactification of $\C^3$ by adjunction of the plane at infinity $\Delta_{\infty}$, and let $M$ stands for the manifold
obtained from $\C \p (3)$ through a punctual blow-up at the origin. The manifold $M$ can be viewed as a fiber bundle by
projective lines equipped with two natural projections, namely
\begin{align*}
	\mathcal{P}_0&: \, \widetilde{\C}^3 \to \Delta_0 \\
	\mathcal{P}_{\infty}&: \, \widetilde{\C}^3 \to \Delta_{\infty} \, ,
\end{align*}
where $\widetilde{\C}^3$ stands for the blow-up of $\C^3$ at the origin and $\Delta_0$ represents the divisor obtained by
the punctual blow-up of $\C \p (3)$ at the origin.

\begin{figure}[hbtp]
	\centering
	\includegraphics[scale=0.75]{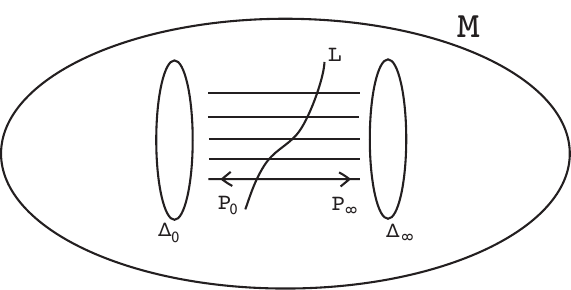}
	\caption{The manifold $M$ and the corresponding bundle projections}
	\label{fig:bundleprojection}
\end{figure}

Since the vector field $X$ is polynomial, it admits a meromorphic extension to $M$, where $\Delta_0$ corresponds to the zero
divisor and $\Delta_{\infty}$ to the pole divisor. The vector field $X$ induces a holomorphic foliation $\tilf$ on $M$ and,
since we are assuming $X$ not to be a multiple of the Radial vector field, $\tilf$ leaves the two divisors $\Delta_0$ and
$\Delta_{\infty}$ invariant. Note that, since $X$ is homogeneous, the projection of every leaf $L$ onto $\Delta_0$ (resp.
$\Delta_{\infty}$), $\calp_0 (L) =L_0$ (resp. $\calp_{\infty} (L) =L_{\infty}$), is clearly a leaf of $\tilf_0$ (resp.
$\tilf_{\infty}$), the restriction of $\fol$ to $\Delta_0$ (resp. $\Delta_{\infty}$). Let then $L$ be non-algebraic leaf
of $\tilf$ not contained in $\Delta_0 \cup \Delta_{\infty}$. We have that the restriction of $\calp_0$ (resp. $\calp_{\infty}$)
to $L$ realizes $L$ as an Abelian covering of $L_0$ (resp. $L_{\infty}$). It then follows that the non-compact leaves $L, \,
L_0, \, L_{\infty}$ have all the same nature: either they are all covered by $\C$ or they are all covered by the unit disc~$D$.
Furthermore $L_0, \, L_{\infty}$ are isomorphic as Riemann surfaces while $L$ is an Abelian covering of $L_0, \, L_{\infty}$.

In order to study the behavior of the solutions nearby the infinity (i.e. away from compact subsets of $\C^3$), let $M$ be
equipped with affine coordinates $(x , y, z)$ such that
\begin{enumerate}
	
	\item[(i)] $\{ z=0\} \subset \Delta_{\infty}$, $(x , y) \in \C^2$, $z \in \C$.
	
	\item[(ii)] the transformed $\tXX$ of the vector field $X$ on $M$ is given by
	\begin{equation}
		\tXX = \frac{1}{z^{d-1}} \left[ F (x , y) \frac{\partial}{\partial
			x} + G (x , y) \frac{\partial}{\partial y} + zH (x , y)
		\frac{\partial}{\partial z} \right]   \label{tXX}
	\end{equation}
	where $F ,G$ are polynomials of degree~either $d$ or $d+1$ and $H$ is a polynomial of degree~$d$ (the independence of $F ,G$
	and $H$ on the variable $z$ is a consequence of the homogeneous character of $X$).

	\item[(iii)] The projection $\calp_{\infty} : M \rightarrow \Delta_{\infty}$ in the above coordinates becomes $(x , y, z)
	\mapsto (x , y)$.
\end{enumerate}
Also, it can be assumed that the line at infinity $\Delta_{\infty}^{(x,y)}$ over the plane at infinity $\Delta_{\infty}
\simeq \{z=0\}$ is not invariant by $\tilf$.

Let $L_{\infty} \subseteq \Delta_{\infty}$ be a leaf of the foliation $\tilf_{\infty}$. Denote by $S$ the cone over $L_{\infty}$,
i.e. let $S = \mathcal{P}_{\infty}^{-1}(L_{\infty})$. The mentioned cone is a $2$-dimensional immersed singular surface clearly
invariant by $\fol$. In particular, $S$ is naturally equipped with a holomorphic foliation, denoted by $\tilf_S$, having a transversely
conformal structure.

\begin{figure}[hbtp]
	\centering
	\includegraphics[scale=0.65]{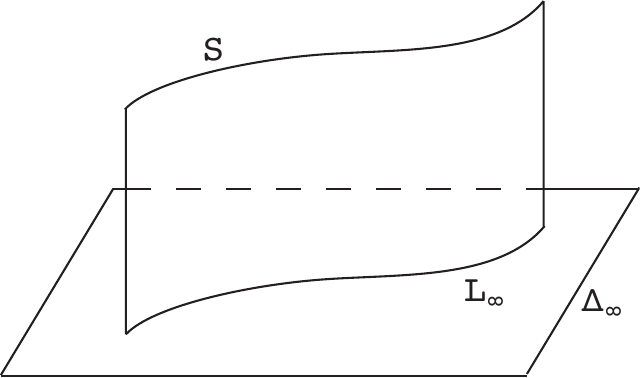}
	\caption{The invariant cone over $L_{\infty}$}
	\label{fig:invariancone}
\end{figure}

Let then $L$ be a leaf of $\fol$ contained in the cone over $L_{\infty}$. The first step of our method consists on having
quantitative estimates on the ``speed'' with which $L$ approaches $\Delta_{\infty}$. As mentioned before, since $\Delta_{\infty} \subseteq
\C\p(n)$ and the Fubini-Study metric on $\C\p(n)$ has positive curvature, it is well known that complex submanifolds always bend
towards $\Delta_{\infty}$ (see for example~\cite{langevin}). So, to keep a ``good control'' of the directions over which the leaf
$L$ approach the infinity we proceed as follows. We shall equip $L$ with an Abelian $1$-form $\omega_1$ related to the holonomy
of the leaf in question. To be more precise, suppose that $L_{\infty}$ is locally parameterized by
\[
x \to (x,y(x),0) \, .
\]
Then $L$ is parameterized by $x \to (x,y(x),z(x))$, where $z = z(x)$ satisfies
\[
\frac{dz}{dx} = \frac{H(x,y(x))}{F(x,y(x))} dx \, .
\]
It then follows that
\[
z = z_0 {\rm exp} \left( \int_{x_0}^x \frac{H(x,y(x))}{F(x,y(x))} dx \right) \, .
\]
We have then that the Abelian $1$-form
\[
\omega_1 = \frac{H(x,y(x))}{F(x,y(x))} dx
\]
is the logarithmic derivative of the holonomy in the sense that if $c$ is a path on $L_{\infty}$ then
\[
\left( {\rm Hol} (c)\right)' (c(0)) = e^{-\int_c \omega_1} \, .
\]

Fix the a point $p \in L_{\infty}$. We claim that there exist real trajectories on $L_{\infty}$ having contractive holonomy.
The trajectories in question correspond to the leaves of the real oriented foliation on $L_{\infty}$ defined by
\[
\mathcal{H} \, : \, \, \, \, \{{\rm Im} (\omega_1) = 0\} \, ,
\]
where ${\rm Im} (\omega_1)$ stands for the imaginary part of $\omega_1$ and the orientation is such that
\[
{\rm Re} (\omega_1 . \phi'(t)) = \omega_1 . \phi'(t) > 0 \, .
\]
In fact, if $c : [0,1] \to l$ is the parametrization of a leaf $l$ of $\mathcal{H}$, $l \subseteq L_{\infty}$, then
\[
\left| {\rm Hol} (c)'\right| = e^{-\int_c \omega_1} = e^{-{\rm Re} \left(\int_c \omega_1\right)} < 1 \, .
\]
Note that the trajectories defined above are not the only trajectories having a contractive holonomy. Is fact, for every
fixed $\theta \in ]-\pi/2, \pi/2[$, the oriented real foliation $\mathcal{H}_{\theta}$ making an angle $\theta$ with
$\mathcal{H}$ is such that the holonomy with respect to their leaves is contractive.

Fix then a leaf $L$ of $\fol$ contained in the cone over $L_{\infty}$. Fix a point $p \in L_{\infty}$ and let $q$ be a
point projecting on $p$. Let $l_p \subseteq L_{\infty}$ be a leaf of $\mathcal{H}$ and let $l_q$ stands for the lift of
the mentioned leaf to $L$. A first remark that can be made is that our lift does not leave the affine coordinates $(x,y,z)$
above. In fact, it can be checked that points in the line at infinity $\Delta_{\infty}^{(x,y)}$ of the plane at infinity
$\Delta_{\infty}$ provide singularities for $\mathcal{H}$ of source type. To prove Theorem~\ref{teo_area} we have to
control the ``high'' (i.e. the distance of $L$ to $L_{\infty}$) and the time that we pass away from a fixed neighborhood
$V$ of the singular set of $\fol$.

\begin{figure}[hbtp]
	\centering
	\includegraphics[scale=0.65]{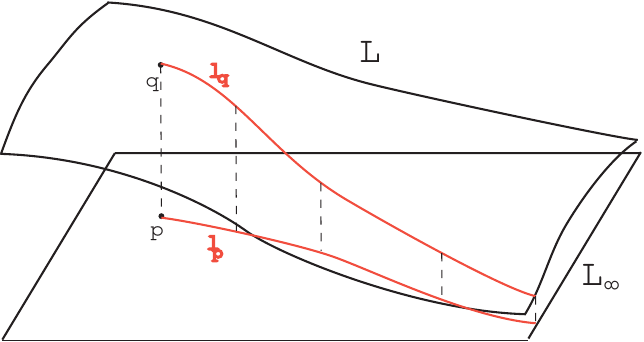}
	\caption{A leaf $L$, its projection $L_{\infty}$ and the leaves of $\mathcal{H}$ passing through the points $p$ and $q$.}
	\label{fig:leaves_of_H}
\end{figure}

To begin with note that away from $V$ the $1$-form $\omega_1$ is bounded from below  by a positive constant $\alpha$ (up
to considering the parametrization $y \to (x(y),y,0)$ instead of the considered one). Note that although we are away from
a fixed neighborhood of the singular set of $\fol$, the domain of definition may contain singularities of the real foliation
$\mathcal{H}$. Singularities of $\mathcal{H}$ may be of three types: sinks, sources or saddles.  The two first ones provide
a minimum or a maximum to the (local) distance from $L_{\infty}$ to $L$, respectively. If a minimum is attained, then we
have ``arrived'' to $\Delta_{\infty}$ since, as mentioned before, we cannot have a minimum unless it is zero. Furthermore,
it is clear that sources are not reached by the (``positive'' direction) of our oriented leaves. Finally, concerning sources,
it can be proved that we can exclude an arbitrarily small neighborhood of it and still keep the contractive holonomy by following
the leaves of $\mathcal{H}_{\theta}$ for some $\theta$ belonging to $]-\pi/1 + \delta, \pi/2 - \delta[$, with $\delta > 0$.

To finish the idea of the proof of Theorem-\ref{teo_area}, let us just show some estimates. To be brief we will explain how to proceed
in the case we stay away from a fixed neighborhood $V$ of the singular set of the foliation. The main idea is to prove that the time passed
in $M \setminus V$ is finite. We recall that the singularities of $\calh$ at points in $\Delta_{\infty}^{(x,y)}$ (the line at infinity
of the hyperplane at infinity) are ``source-like'' so that an oriented trajectory of $\calh$ cannot actually intersect $\Delta_{\infty}^{(x,y)}$.
Though these trajectories of $\calh$ may come ``close'' to $\Delta_{\infty}^{(x,y)}$, owing to Lemma~3.10 of \cite{RR_applications} we
know that every sufficiently long segment of $l_p$ has ``most of its length'' contained in a fixed compact part of the affine $\C^2$
associated to the coordinates $(x,y)$. Let then a compact part $K$ of the mentioned affine copy of $\C^2$ be fixed and let us precise
the estimates we need on this compact part - the estimates of Lemma~3.10 about the non-compact part allows us to adapt the estimate we
present below away from $K$.

So, let $c:[0,1] \to l_q$ be a parametrization of a connected path of $l_q$ above $l_p$. We have that the ``high'' $z = z(t)$ along $l_q$
satisfies
\begin{align*}
	|z| &= \left| z_0 e^{-\int_c \omega_1} \right| = |z_0| e^{-{\rm Re} (\int_c \omega_1)} \\
	&= |z_0| e^{-\int_0^1 {\rm Re} (\omega_1(c(t)). c'(t)) dt} = |z_0| e^{-\int_0^1 |\omega_1(c(t)). c'(t)| dt} \\
	&\leq |z_0| e^{-\alpha . {\rm lenght} (c)} \, ,
\end{align*}
where the last inequality comes from the fact that $\omega_1$ is bounded from below by $\alpha$ away from $V$. We have then
that if the length of $c$ is greater than $\ln 2 / \alpha$, then
\[
|z_1| = |z(1)| \leq \frac{|z(0)|}{2} = \frac{|z_0|}{2} \, .
\]

We have finally to control the time we take to cover the path $l_q$. Recall that the time-form associated to a leaf $L$ is
well-defined provided that $L$ is not contained in the divisor of zeros and poles of $X$. If the vector field $X$ is supposed
to be semicomplete, then its restriction to $L$ is everywhere holomorphic and the orders of its zeros cannot exceed~$2$. It
follows at once that $dT$ is meromorphic on all of $L$ and it has no zeros. Furthermore, the poles of $dT$ have order bounded
by~$2$. Finally, recall also that given a curve $c : [0,1] \to L$ joining two points $c(0)$ and $c(1)$ in $L$ satisfying $X(c(0))
\ne 0$ and $X(c(1)) \ne 0$, the integral $\int_c dT$ measures the time needed to traverse $c$ from $c(0)$ to $c(1)$ following the
flow of $X$ as long as $X$ is semicomplete. In fact, when a vector field is semicomplete the notion of time arising from its
semi-global flow is well-defined.

Thus the integral $\int_{l_q} dT$ can be estimated as follows. The time-form on $L$ is given, in local coordinates $(x,y,z)$, by
$dT = z^{d-1}dx/F(x,y)$. Since $l_p$, the projection of $l_q$ by $\mathcal{P}_{\infty}$, is contained on a compact set
not intersecting the singular set of $\tilf_{\infty}$, the absolute value of $F(x,y)$ is bounded from below, i.e. $|F(x,y)|
\geq \db > 0$ for all$(x, y) \in \Delta_{\infty} \setminus V$ and some positive constant $\beta$. Otherwise we replace $F$ by $G$
(recall that we are dealing only with regular points of $\tilf$ on $\Delta_{\infty}$). Hence, considering $l_q$ as the
concatenation of segments having length equal to $\ln(2)/\da$, $l_q = \sum_{i=0}^{\infty}l_{i,q}$, it follows that
\begin{eqnarray*}
	\left| \int_{l_q} dT \right| &=& \left| \sum_{i=0}^{\infty} \int_{l_{i,q}} \frac{z^{d-1}}{F(x, y)} dx \right|
	\leq  \sum_{i=0}^{\infty} \left| \int_0^1 \frac{z_{i,q}^{d-1}(t)}{F(x_{i,q}(t),y_{i,q}(t))} x^{\prime}_{i,q}(t) dt \right|\\
	&\leq&  \sum_{i=0}^{\infty} \int_0^1 \frac{|z_{i,q}(t)|^{d-1}}{|F(x_{i,q}(t),y_{i,q}(t))|} |x^{\prime}_{i,q}(t)| dt
	\leq  \sum_{i=0}^{\infty} \int_0^1 \frac{|z_0|^{d-1}(\frac{1}{2})^{i(d-1)}}{\db} |l^{\prime}_{i,p}(t)| dt\\
	&\leq&  \frac{|z_0|^{d-1}}{\db} {\rm length }(l_{i,p}) \sum_{i=0}^{\infty}\left( \frac{1}{2^{d-1}} \right)^{i}
	=  \frac{|z_0|^{d-1} {\rm ln}(2)}{\da \db} \frac{1}{1 - (\frac{1}{2})^{d-1}} < \infty
\end{eqnarray*}
where $l_{i,q}(t) = (x_{i,q}(t), y_{i,q}(t), z_{i,q}(t))$, $t \in [0,1]$, is such that $l_q = \sum_{i=0}^{\infty} l_{i,q}$ and
$\mathcal{P}_{\infty}(l_{i,q}) = l_{i,p}$. The conclusion follows.


\section{Vector fields with univalued solutions - Local aspects}\label{sec:localaspects}

The interest of the notion of semicomplete vector fields (cf. Section~\ref{sec:sc_global_dynamics}) comes from the fact that
the restriction of a complete holomorphic vector field, defined on a manifold $M$, to every open set $U \subseteq M$, is
automatically semicomplete. Furthermore, given a semicomplete vector field on an open set $U$, its restriction to an open
set $V \subseteq U$ is also semicomplete. In this sense, we are allowed to talk about germs of semicomplete vector fields.
With abuse of notation, we can also talk about semicomplete singularities. By semicomplete singularities we mean a singular
point associated to a semicomplete vector field on a small neighborhood of the point in question. Note however that a singularity
may have more than one representative vector field and not all of them need to be semicomplete.

From the preceding it follows that semicomplete vector fields can be viewed as the ``local version'' of complete vector fields.
In fact, a singularity that is not semicomplete cannot be realized by a complete vector field. In particular, it cannot be
realized by a globally defined holomorphic vector field on a compact manifold. The understanding of semicomplete singularities
is then important to the understanding of holomorphic vector fields (globally) defined on compact manifolds.

There is a long standing well-known question by E. Ghys that can be formulated in terms of semicomplete vector fields as follows:

\begin{quest}[{\bf Ghys' question}]
	Let $X$ be a holomorphic vector field on $(\C^n,0)$ that is semicomplete and having an isolated singularity at the origin. Is it true that
	$J^2 X (0) \ne 0$, i.e. must the second jet of $X$ at the singular point be different from zero?
\end{quest}

His motivation is, at least partially, related to problems about bounds for the dimension of automorphism group of compact complex manifolds. To be more precise,
consider a compact complex manifold $M$ and denote by ${\rm Aut} \, (M)$ the group of holomorphic diffeomorphisms
of $M$. It is well-known that ${\rm Aut} \, (M)$ is a finite dimensional complex Lie group whose Lie algebra can be identified
with $\mathfrak{X} \, (M)$, the space of all holomorphic vector fields defined on $M$. It is also known that the dimension of the automorphism
group of $M$ cannot be bounded by the dimension of $M$, in general: it is sufficient to consider the family of Hirzebruch surfaces
$\{F_n\}$, whose dimension of automorphism group equals $n+5$, for $n \geq 1$. However the same question can be formulated for
special classes of manifolds. For example, there is a question attributed to Hwuang and Mok that asks if $\C\p(n)$ is the projective manifold with the largest (in terms of dimension) automorphism group among manifolds with Picard group isomorphic to $\Z$.

Let us briefly explain how an affirmative answer to Ghys conjecture can help us in the above mentioned problems. Suppose that $M$ is
a compact complex manifold of dimension $n$ and fix a point $p \in M$ and $k \in \N$. We have the following short exact sequence
\[
\mathfrak{X}_p^k \, (M) \to \mathfrak{X} \, (M) \to J^k_p(M) \, ,
\]
where $\mathfrak{X}_p^k(M)$ stands for the set of holomorphic vector fields with vanishing $k$-jet at $p$ and $J^k_p(M)$ denotes
the space of $k$-jets. Thus, we have
\[
{\rm dim} \, \mathfrak{X} \, (M) \leq {\rm dim} \, \mathfrak{X}_p^k \, (M) + {\rm dim} \, J^k_p(M) \, .
\]
Concerning the space of jets, there exist effective bounds for ${\rm dim} \, J^k_p(M)$ in terms of $n = {\rm dim} \, (M)$. If
${\rm dim} \, \mathfrak{X}_p^k \, (M)$ has bounds in terms of ${\rm dim} \, (M)$ for a certain $p \in M$ and $k \in M$, then
${\rm dim} \, \mathfrak{X} \, (M)$ and, consequently, ${\rm dim} \, ({\rm Aut} \, (M))$ has such bounds as well. For example,
suppose that we happen to know that for a certain class of manifolds every singularity of a vector field is necessarily
isolated. Then, if Ghys conjecture hold then ${\rm dim} \, \mathfrak{X}_p^3 \, (M) = 0$ and then we obtain the desired bounds
for ${\rm dim} \, \mathfrak{X}$.

So, let us focus on Ghys conjecture. In the paper~\cite{Rebelo96}, where the notion of semicomplete vector field was introduced, the following has been proved.

\begin{teo}\cite{Rebelo96}
	Let $X$ be a vector field defined on a neighborhood of the origin of $\C^2$ with an isolated singularity at the origin. If $X$ is semicomplete, then $J_0^2 X \ne 0$.
\end{teo}

The proof relies on the fact that the foliation associated to $X$ possesses at least one separatrix, i.e. a germ of analytic curve
that is invariant by the foliation in question. Being the singular set of $X$ reduced to the origin, it follows that the restriction of $X$ to the mentioned
invariant curve does not vanish identically. Furthermore, this restriction is still a semicomplete vector field. Considering then the restriction of $X$ to
a separatrix the problem is essentially reduced to an one-dimensional situation. Note that, in the case where the separatrix is singular, a Puiseaux parameterization for the separatrix can be used to settle this issue. This one-dimensional case
is treated in the same paper by direct methods.

Later, semicomplete vector fields with an isolated singular point at the origin and vanishing linear part have been characterized (cf.~\cite{RG}). A
geometric study of the above mentioned vector fields/corresponding foliations has been done allowing the authors to prove that all those vector fields
are integrable in the sense that they admit a (non-trivial) holomorphic or meromorphic first integral. Furthermore, the vector field $X$ is conjugate
to its first non-zero homogeneous component thus providing a sharp classification theorem. Recall that this has already been mentioned in Section~\ref{sec:sc_global_dynamics} when Theorem~\ref{teo:sc_quadratic} was stated.

The question of whether or not the above results still hold for semicomplete vector fields (at isolated singular points) in higher dimensional
manifolds is a natural one. However, it is easy to detect a number of new difficulties that will play a role in any attempt at generalizing the preceding results to higher dimensional manifolds. Let us enumerate some of them.
\begin{itemize}
	\item[(1)] Unlike the case of holomorphic foliations on $(\C^2,0)$, there exist
	holomorphic foliations on $(\C^3,0)$ with an isolated singular point but no separatrix (i.e. invariant analytic curve through the singular point,
	cf. Section~\ref{sec:invariant_sets}). To prove, for example, that these vector fields without separatrices fail to be semicomplete is a challenging problem.
	
	\item[(2)] Another ingredient that played an important role on the classification of semicomplete holomorphic vector field in dimension~$2$ is the resolution
	theorem of Seidenberg. The lack of a faithfully analogous procedure for reducing the singularities
	of vector fields in dimensions~$3$ and higher (cf. Section~\ref{Sec:resolution} for recent results)
	adds therefore to the difficulty of the problem.
	
	\item[(3)] Even in the case were we are given a holomorphic foliation admitting a simple reduction of singularities in
	the sense of Seidenberg (where the linear part of the blown-up foliation has at least a non-vanishing eigenvalue at each singular point)
	additional difficulties are expected if we compare with the two-dimensional case. In fact, already in the three-dimensional case saddle-nodes
	of codimension~$2$ (i.e. with two eigenvalues equal to zero) may appear and these singularities are still poorly understood.
\end{itemize}

The first deep investigations involving semicomplete vector fields in higher dimensions were conducted by
A. Guillot in~\cite{guillotFourier}, and \cite{guillotIHES}. These investigations soon confirmed that the case of, say dimension~$3$, was
far more subtle than its two-dimensional version.
Indeed, in the mentioned papers by A. Guillot there are a
huge variety of examples of semicomplete vector fields with isolated singularities and exhibiting interesting dynamical properties. As to genuinely complicated dynamical
behavior, A. Guillot obtained some remarkable examples by studying Halphen vector
fields. An exhaustive
classification of all semicomplete vector fields
analogous to the classification given in~\cite{RG} is therefore unlikely or, at least, not particularly useful. Furthermore, from the characterization of their dynamics, it can be concluded that they do not admit holomorphic/meromorphic first integrals.
Guillot's work was extended to the meromorphic setting allowing for additional dynamical complications in a joint work with A. Elshafei and J. Rebelo (see \cite{Ahmed_thesis} and \cite{Ahmed_paper} for details). Here we should recall that a meromorphic vector field is said to be semicomplete if it is semicomplete away from its pole divisor.

Whereas the
papers~\cite{Rebelo96} and~\cite{RG} deal with isolated singularities of a $\C$-action or, equivalently, of a complete vector field
$X$ on a complex two-dimensional ambient space, there was significant evidence that a ``natural'' extension of the methods and results obtained in these papers
in dimension~$2$ might be achieved by considering two commuting vector fields or, more precisely, $\C^2$-actions of rank~$2$ (again the reader will note
that a singularity of a globally defined $\C^2$-action is automatically semicomplete).

Let us focus on the problem about the vanishing order of a semicomplete vector field at an isolated singular point. The general principle
to Ghys conjecture is the existence of a separatrix through the isolated singular point.

\begin{prop}
	Let $X$ be a semicomplete vector field defined on a neighborhood of the origin of $\C^n$ and having an isolated singularity at the origin. If $X$ admits a separatrix through the origin, then $J_0^2 X \ne 0$.
\end{prop}

Recall that a holomorphic
foliation by curves on a complex surface always admits separatrices through its singular points but this no longer holds when the ambient
manifold is of dimension~$3$ or greater. However, as mentioned in Section~\ref{sec:invariant_sets}, it has been proved in~\cite{RR_secondjet} that in the case we are given
two holomorphic vector fields immersed in a representation of a Lie algebra of dimension~2, that are in addition linearly independent up to
a set of codimension at least~$2$, then $X, \, Y$ possess a common separatrix. As an important consequence of this result we were able to
prove in the same paper that Ghys conjecture holds for vector fields on $3$-dimensional compact manifolds whose automorphism group has
dimension at least~$2$. More precisely, the following has been proved:

\begin{teo}\cite{RR_secondjet}\label{teo_RR_nonvanishing_secondjet}
	Consider a compact complex manifold $M$ of dimension~$3$ and assume that the dimension of ${\rm Aut}\, (M)$ is at least~$2$. Let $Z$ be an
	element of $\mathfrak{X}\, (M)$ and suppose that $p \in M$ is an isolated singularity of~$Z$. Then
	$$
	J^2 (Z) \, (p) \neq 0
	$$
	i.e., the second jet of $Z$ at the point $p$ does not vanish.
\end{teo}

\begin{rem}
	{\rm Note that $M$ is not assumed to be algebraic in the above theorem. We should mention that in the case where $M$ is algebraic, or more generally K\"ahler, then
		the statement holds in arbitrary dimensions and regardless of the condition on the dimension of ${\rm Aut}\, (M)$.}
\end{rem}

Our technique to derive Theorem~\ref{teo_RR_nonvanishing_secondjet} from Theorem~\ref{teo_RR_secondjet} also yields another interesting
result. In fact, the assumption that $M$ is compact is not fully indispensable in many cases. For example, suppose that $N$ is a
Stein manifold of dimension~$3$ and suppose that $N$ is effectively acted upon by a finite dimensional Lie group $G$. Then the Lie algebra
$\mathfrak{G}$ of $G$ embeds into the space $\mathfrak{X}_{\rm comp} (N)$ of {\it complete}\, holomorphic vector fields on $N$. The
study of these complete holomorphic vector fields on Stein manifolds is a topic of interest having its roots in a classical
work of Suzuki~\cite{suzuki}. In this direction, our techniques yield:

\begin{teo}
	Let $N$ denote a Stein manifold of dimension~$3$ and consider a finite dimensional
	Lie algebra $\mathfrak{G}$ embedded in $\mathfrak{X}_{\rm comp} (N)$ (the space of complete holomorphic vector fields on $N$).
	Assume that the dimension of $\mathfrak{G}$ is at least~$2$. If $Z$ is an element of $\mathfrak{G} \subseteq
	\mathfrak{X}\, (M)$ possessing an isolated singular point $p \in N$, then the linear part of $Z$ at $p$ cannot vanish,
	i.e. $p$ is a non-degenerate singularity of~$Z$.
\end{teo}


\section{Generic pseudogroups on \texorpdfstring{$(\C ,0)$}{c0}  and the topology of leaves}\label{Sec_topology}

In the study of some well-known problems about singular holomorphic foliations, we usually experience difficulties
concerning to greater or lesser extent the topology of their leaves. Yet, most of these problems are essentially
concerned with pseudogroups generated by certain local holomorphic diffeomorphisms defined on a neighborhood of $0
\in \C$ (recall definition below). In this sense, results about pseudogroups of $\diff$ generated by a finite number of
local holomorphic diffeomorphisms are crucial for the understanding of certain singular foliations defined about the
origin of $\C^2$. Furthermore, for most of these problems, it is necessary to consider classes of pseudogroups with
a distinguished generating set all of whose elements have fixed conjugacy class in $\diff$. These statements will be
explained below using a standard example.

Contrary to the previous sections, we consider here a singular holomorphic foliation defined about the origin in $\C^2$
and recall that these foliations are obtained by means of holomorphic vector fields having isolated singular points. The
study of these singularities and of their deformations, paralleling Zariski problem, led to the introduction of the
{\it Krull topology}\, in the space of these foliations. In this topology, a sequence of foliations $\fol_i$ is said
to converge to $\fol$ if there are representatives $X_i$ for $\fol_i$ and $X$ for $\fol$ such that $X_i$ is tangent to
$X$, at the origin, to arbitrarily high orders (modulo choosing $i$ large enough). It should be noted that, given a
foliation $\fol$, its resolution depends only on a finite jet of the Taylor series of $X$ at the singular point in the
sense that if $\fol'$ is close to $\fol$ in the Krull topology, both foliations will admit the same resolution map. Furthermore,
the position of the singularities of the resolved foliations $\tilf, \, \tilf'$ will also coincide as well as the corresponding
eigenvalues.

As an example, consider a nilpotent foliation $\fol$ associated to Arnold singularity $A^{2n+1}$, i.e. a nilpotent foliation
admitting a unique separatrix that happens to be a curve analytically equivalent to the cusp of equation $\{ y^2 - x^{2n+1} = 0\}$.
An important remark from what precedes is that whenever $\fol'$ is sufficiently close to $\fol$ in the Krull topology, $\fol'$ is
also a nilpotent foliation of type $A^{2n+1}$. In other words, the class of Arnold singularities is closed for small perturbations in
the Krull topology. It is then natural to wonder what type of dynamical behavior can be expected from these foliations, or more
precisely, from a ``typical'' foliation in this family. The following are examples of long-standing problems in the area:
\begin{itemize}
	\item[(1)] Does there exist a nilpotent foliation $\fol$ in $A^{2n+1}$ whose leaves are simply connected (apart maybe from a countable set)?
	
	\item[(2)] Is the set of these foliations dense in the Krull topology, i.e. given a nilpotent foliation $\fol$ in $A^{2n+1}$,
	does there exist a sequence of foliations $\fol_i$ converging to $\fol$ in the Krull topology and such that every $\fol_i$ has
	simply connected leaves (with possible exception of a countable set of leaves)?
\end{itemize}

Our methods in \cite{MRR} are powerful enough to affirmatively settle both questions above. Moreover, in the paper~\cite{RR_stabilizers}
we also establish that the countable set is, indeed, infinite and that the non-simply connected leaves are topologically cylinders.
More precisely, in these two papers it is proved the following:

\begin{teo}\cite{MRR,RR_stabilizers}\label{foliation}.
	Let $X \in \mathfrak{X}_{(\C^2,0)}$ be a vector field with an isolated singularity at the origin and defining
	a germ of nilpotent foliation $\fol$ of type $A^{2n+1}$. Then, for every $N \in \N$, there exists a vector field
	$X' \in \mathfrak{X}_{(\C^2,0)}$ defining a germ of foliation $\fol^{\prime}$ and satisfying the
	following conditions:
	\begin{itemize}
		\item[(a)] $J_0^N X^{\prime} = J_0^N X$.
		
		\item[(b)] $\fol$ and $\fol^{\prime}$ have $S$ as a common separatrix.
		
		\item[(c)] there exists a fundamental system of open neighborhoods $\{U_j\}_{j \in \N}$ of $S$, inside a
		closed ball $\bar{B}(0,R)$, such that the following holds for every $j \in \N$:
		\begin{itemize}
			\item[(c1)] The leaves of the restriction of $\fol^{\prime}$ to $U_j \setminus S$, $\fol^{\prime}|_{(U_j \setminus S)}$
			are simply connected except for a countable number of them.
			
			\item[(c2)] The countable set constituted by non-simply connected leaves is, indeed, infinite.

			\item[(c3)] Every leaf of $\fol^{\prime}|_{(U_j \setminus S)}$ is either simply connected or homeomorphic
			to a cylinder.
		\end{itemize}
	\end{itemize}
\end{teo}

These problems are related with pseudogroups which is generated by certain elements of  $\diff$, as it will be explained in the next paragraph. In this sense, let us start by recalling the notion of pseudogroup. Consider the group $\diff$ of germs of
holomorphic diffeomorphisms fixing $0 \in \C$, where the group law is induced by composition. Assume that $G$ is actually
a subgroup of $\diff$ generated by the elements $h_1, \ldots, h_k$. Then, consider a small neighborhood $V$ of the origin
where the local diffeomorphisms $h_1, \ldots, h_k$, along with their inverses $h_1^{-1}, \ldots, h_k^{-1}$, are all well
defined diffeomorphisms onto their images. The {\it pseudogroup} generated by $h_1, \ldots, h_k$ (or rather by $h_1, \ldots,
h_k, h_1^{-1}, \ldots, h_k^{-1}$ if there is any risk of confusion) on $V$ is defined as follows. Every element of this
pseudogroup has the form $F = F_s \circ \ldots \circ F_1$ where each $F_i$, $i \in \{1, \ldots, s\}$, belongs to the set
$\{h_i^{\pm 1}, i=1, \ldots, k\}$. The element $F$ should be regarded as an one-to-one holomorphic map defined on a subset
of $V$. Indeed, the domain of definition of $F = F_s \circ \ldots \circ F_1$, as an element of the pseudogroup, consists
of those points $x \in V$ such that for every $1 \leq l < s$ the point $F_l \circ \ldots \circ F_1(x)$ belongs to $V$.
Since the origin is fixed by the diffeomorphisms $h_1, \ldots, h_k$, it follows that every element $F$ in this pseudogroup
possesses a non-empty open domain of definition. This domain of definition may however be disconnected. Whenever no
misunderstanding is possible, the pseudogroup defined above will also be denoted by $G$ and we are allowed to shift
back and forward from $G$ viewed as pseudogroup or as group of germs.

\begin{figure}[hbtp]
	\centering
	\includegraphics[scale=0.75]{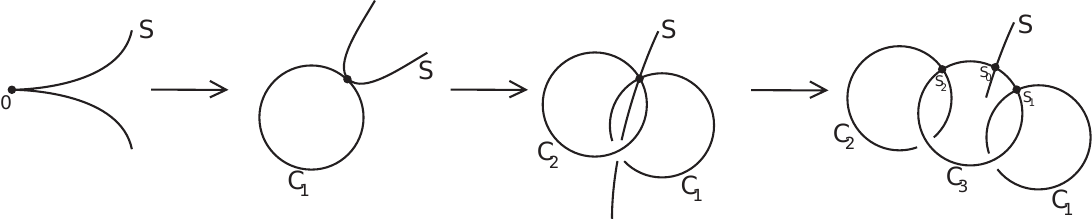}
	\caption{The desingularization diagram of the foliation associated with $A^3$}
	\label{Resolucao}
\end{figure}

So, let us explain how the above problems are concerned with pseudogroups generated by certain local elements of $\diff$.
To do this we will restrict ourselves to the particular case of a nilpotent foliation $\fol$ associated to Arnold singularity
$A^3$, i.e. a nilpotent foliation admitting a unique separatrix $S$ that happens to be a curve analytically equivalent to the
cusp of equation $\{ y^2 - x^3 = 0\}$. For this type of foliation, the map associated to the desingularization of the separatrix
$E_S: M \rightarrow \C^2$ reduces also the foliation $\fol$ (see Figure~\ref{Resolucao} for the corresponding resolution). So, let us
describe the resolution in question.

To begin with, let us consider standard coordinates $(x,y)$ for $(\C^,0)$ where the separatrix if given by $\{y^2 - x^3 = 0\}$.
The origin of the mentioned coordinates is the (unique) singular point of $\fol$. So, let us first consider the one-point blow-up
of $\fol$ centered at the origin. The strict transform of the separatrix is tangent to the resulting component of the exceptional
divisor, denoted by $C_1$, at some point. This point of tangency is the unique singular point for the strict transform foliation
and it is a degenerate singular point. So, consider now the punctual blow-up of the transformed foliation at this new singular point.
Let $C_2$ be the resulting irreducible component of the exceptional divisor. Since $C_1$ and $S$ were tangent, it follows that $C_1$,
$C_2$ and $S$ intersect all each other at the same point. This intersection is however transverse at this point. Nonetheless the
eigenvalues of the foliation at this intersection and singular point are both equal to zero and so we need to perform a punctual
blow-up at this point. Let us then perform a third one-point blow-up, centered at this intersection point and let $C_3$ be the
new irreducible component of the exceptional divisor. Now, since $C_1, C_2$ and $S$ were transverse we have that their strict
transforms intersect $C_3$ at distinct points. We denote by $s_1$ the intersection point of $C_1$ with $C_3$, by $s_2$ the
intersection point of $C_2$ with $C_3$ and by $s_0$ the intersection point of the separatrix with $C_3$.

All those singular points $s_0, s_1$ and $s_2$ are non-degenerate singular points for $\fol$. In fact, we have that both eigenvalues
of the transformed foliation, still denoted by $\fol$, at each one of those singular points are different from zero. For example, it
can easily be checked that
\begin{itemize}
	\item[1.] the eigenvalues of $\fol$ at $s_1$ are $1, -3$;
	\item[2.] the eigenvalues of $\fol$ at $s_2$ are $1, -2$.
\end{itemize}
The eigenvalues of $\fol$ at $s_0$ can also be deduced by using the index formula.

Clearly, every component of the exceptional divisor is invariant by the foliation in question. So, $C_1 \setminus \{s_1\}$
is a leaf of $\fol$ and so is $C_2 \setminus \{s_2\}$. We have that $C_1 \setminus \{s_1\}$ is isomorphic to $\C$ and thus,
it is simply connected. This means that the holonomy of $\fol$ with respect to this leaf is the identity. Now, recalling that
the quotient of the eigenvalues of $\fol$ at $s_1$ is negative real, it follows from Mattei-Moussu that $\fol$ is linearizable
in a sufficiently small neighborhood of $s_1$. Applying the same argument to $C_2 \setminus \{s_2\}$, we get that $\fol$ is
also linearizable in a small neighborhood of $s_2$ as well.

Denoting by $h_{\sigma_1}$ and by $h_{\sigma_2}$ the holonomy map of $\fol$ with respect to small loops on $\C_3$ around $s_1$
and $s_2$, respectively, we have that $h_{\sigma_1}$ is periodic of period~$3$ and $h_{\sigma_2}$ is periodic of period~$2$.
So, $h_{\sigma_1}$ is analytically conjugated to a rotation of angle $2\pi/3$ while $h_{\sigma_2}$ is analytically conjugated to
a rotation of angle $\pi$. So, being $\sigma_0, \sigma_1$ and $\sigma_2$ loops on $C^3$ around $s_0, s_1$ and $s_2$, respectively,
and recalling that they satisfy the relation
\[
\sigma_0 \circ \sigma_1 \circ \sigma_2 = {\rm id} \, ,
\]
it follows that the fundamental group of $C_3$ minus the three singular points is generated by $\sigma_1$ and $\sigma_2$ and,
then, the global holonomy of $\fol$ is generated by $h_{\sigma_1}$ and $h_{\sigma_2}$. In other words, the global holonomy or,
more precisely, the holonomy pseudogroup is generated by the (local) holonomy maps $h_{\sigma_1}$ and $h_{\sigma_2}$ at $s_1$
and $s_2$ w.r.t. the irreducible component intersecting the strict transform of the separatrix. Furthermore, these holonomy
maps are of finite orders ($2$ and~$3$) and hence are linearizable, {\it though not necessarily in the same coordinate}.

From the above paragraph, it follows that every foliation associated to the Arnold singularity $A^3$ gives rise to a pseudogroup
generated by two elements of $\diff$: one having order~$2$ and another having order~$3$. In \cite{MRR}, we proved that the converse
still holds. On other words, we prove that if we are given two diffeomorphisms $f$ and $g$, one being conjugated to a rotation of
order~$2$ and the other one conjugated to a rotation of order~$3$, there exists a foliation as above realizing $f$ and $g$ as
generators of the global holonomy (of holonomy pseudogroup). In this sense the study of the foliations in question is ``equivalent''
to the study of pseudogroups of $\diff$. The proofs of the above theorems are thus reduced to analogous statements for
finitely generated pseudogroups of $\diff$.

\begin{figure}[hbtp]
	\centering
	\includegraphics[scale=0.7]{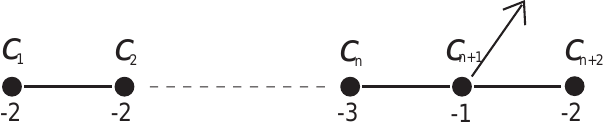}
	\caption{The desingularization diagram of the foliation associated with $A^{2n+1}$}
	\label{graph}
\end{figure}

More generally, for a nilpotent foliation associated to Arnold singularity $A^{2n+1}$, we still have that the holonomy pseudogroup
is generated by two diffeomorphisms $f$ and $g$, one being conjugated to a rotation of order~$3$ but the other one conjugated to a
rotation of order~$2n+1$. In fact, the resolution diagram for such foliation is the same for the corresponding separatrix and is
represented in Figure~\ref{graph}.

The vertices of this graph correspond to the irreducible components of the resulting exceptional divisor. The weight of each irreducible
component equals its self-intersection. In turn, the edges correspond to the intersection of two irreducible components whereas the arrow
corresponds to the intersection point of the (unique) component $C_{n+1}$ of self-intersection~$-1$ with the transform of the separatrix
$S$. The component $C_{n+1}$ contains, as in the previous case, three singular points that we still denote by $s_0, s_1$ and $s_2$ and
where $s_0$ is still the point determined by the intersection of $C_{n+1}$ with the separatrix. Finally $s_1$ (resp. $s_2$) is the
intersection point of $C_{n+1}$ with $C_{n+2}$ (resp. $C_n$).

The singular points of $\fol$ are the intersection points of two consecutive components in the chain $C_1, \ldots, C_{n+2}$ along with
the point $s_0$. All these singular points are simple in the sense that they possess two eigenvalues different from zero. The corresponding
eigenvalues can precisely be determined by using the weights of the various components of the exceptional divisor. The argument applied
in the case of the Arnold singularity $A^3$ allows us to say that the holonomy of $\fol$ associated to the component $C_{n+2}$ (i.e. the
holonomy map associated to the regular leaf $C_{n+2} \setminus \{S_1\}$ of $\fol$) coincides with the identity (the leaf in question is
simply connected). Therefore the germ of $\fol$ at $s_1$ admits a holomorphic first integral and since the corresponding eigenvalues are
$1, 2$, we conclude that the local holonomy map $g$ associated to a small loop around $s_1$ and contained in $C_{n+1}$, has order equal
to~$2$. A similar discussion applies to the component $C_1$ and leads to the conclusion that the local holonomy map $f$ associated to a
small loop around $s_2$ and contained in $C_{n+1}$ has order equal to~$2n + 1$. Since $C_{n+1} \setminus \{s_0, s_1, s_2\}$ is a regular
leaf of $\fol$, we conclude that the (image of the) holonomy representation of the fundamental group of $C_{n+1} \setminus \{s_0, s_1, s_2\}$
in ${\rm Diff} \, (\C, 0)$ is nothing but the group generated by $f,g$. Conversely, given two local diffeomorphism $f, g$ of orders
respectively $2, 2n+1$, they can be realized (up to simultaneous conjugation) as the holonomy of the corresponding component $C_{n+1}$
for some foliation associated to Arnold singularity $A^{2n+1}$. This is done through a well-known gluing procedure explained in~\cite{MRR}.

\begin{figure}[hbtp]
	\centering
	\includegraphics[scale=0.8]{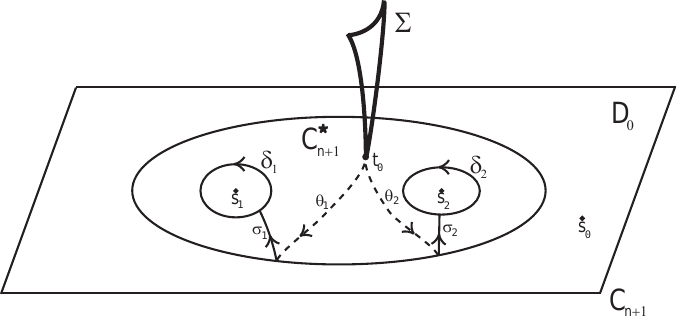}
	\caption{The holonomy representation}
	\label{holonomia}
\end{figure}

Finally, note that the above conclusion depends only on the configuration of the reduction tree which, in turn, is determined by some finite
order jet of~$X$. Hence, if the coefficients of Taylor series of the vector field $X$ are perturbed starting from a sufficiently high order,
the new resulting vector field $X'$ will still give rise to a foliation whose singularity is reduced by the same blow-up map associated to the
divisor of Figure~\ref{graph}. In particular, the holonomy representation of the fundamental group of $C_{n+1} \setminus \{s_0, s_1, s_2\}$
in ${\rm Diff} \, (\C, 0)$, obtained from this new foliation, is still generated by two elements of ${\rm Diff} \, (\C, 0)$ having finite orders
respectively equal to~$2$ and to~$2n+1$. Since every local diffeomorphism of finite order is conjugate to the corresponding rotation, it
follows that the mentioned perturbations are made inside the conjugacy classes of $f$ and $g$. This also justifies the fact that in
Theorems~\ref{teo_diffeo1} and~\ref{teo_diffeo2} below only perturbations of local diffeomorphisms that do not alter the corresponding
conjugation classes were allowed.

So, let $\diff$ be equipped with the so-called analytic topology which, unlike the Krull topology, has the Baire property. Next, consider
a $k$-tuple of local holomorphic diffeomorphisms $f_1, \ldots ,f_k$ fixing $0 \in \C$. Here we impose the condition that the local
diffeomorphisms $f_i$ can be perturbed only inside their conjugacy classes so as to be able to recover results for the initial foliations
in the case where they have finite orders (and hence are all conjugate to a fixed rational rotation).

Fixed $\alpha \in \N$, in the sequel we denote by $\diffalpha$ the subgroup of $\diff$ whose elements are tangent to the
identity to the order~$\alpha$. Finally, we have proved the following.

\begin{teo}\cite{MRR}\label{teo_diffeo1}
	Fixed $\alpha \in \N$, let $f_1, \ldots ,f_k$ be given elements in $\diff$ and consider the cyclic groups $G_1, \ldots , G_k$
	that each of them generates. There exists a $G_{\delta}$-dense set $\mathcal{V} \subset (\diffalpha)^k$ such that, whenever $(h_1, \ldots ,h_k) \in \mathcal{V}$,
	the following holds:
	\begin{enumerate}
		
		\item The group generated by $h_1^{-1}\circ f_1 \circ h_1, \ldots , h_k^{-1}\circ f_k \circ h_k$ induces a group
		in $\diff$ that is isomorphic to the free product $G_1 \ast \cdots \ast G_k$.
		
		\item Let $f_1, \ldots ,f_k$ and $h_1, \ldots ,h_k$ be identified to local diffeomorphisms defined about $0\in
		\C$. Suppose that none of the local diffeomorphisms $f_1, \ldots ,f_k$ has a Cremer point at $0 \in \C$. Denote by
		$\Gamma^h$ the pseudogroup defined on a neighborhood $V$ of $0 \in \C$ by the mappings $h_1^{-1}\circ f_1 \circ
		h_1, \ldots , h_k^{-1}\circ f_k \circ h_k$, where $(h_1, \ldots ,h_k) \in \mathcal{V}$. Then $V$ can be chosen so
		that, for every non-empty reduced word $W (a_1, \ldots , a_k)$, the element of $\Gamma^h$ associated to $W
		(h_1^{-1}\circ f_1 \circ h_1, \ldots , h_k^{-1}\circ f_k \circ h_k)$ does not coincide with the identity on any
		connected component of its domain of definition.
		
	\end{enumerate}
\end{teo}

Theorem~\ref{teo_diffeo1} implies items (a), (b) and (c) of Theorem~\ref{foliation}. It should be noted that perturbing
the foliation inside the class of Arnold singularities of type $A^3$ is equivalent to keeping the conjugacy class of the
generators of their local holonomy maps fixed. Furthermore, the analytic topology allows us to obtain information on
the coefficients of the representative vector fields so as to be able to derive information concerning the Krull topology.
Items (d) and (e) follow from the following result proved in \cite{RR_stabilizers}:

\begin{teo}\cite{RR_stabilizers}\label{teo_diffeo2}
	Suppose we are given $f,g$ in $\diff$ and denote by $D$ an open disc about $0 \in \C$ where $f, g$ and their
	inverses are defined. Assume that none of the local diffeomorphisms $f, g$ has a Cremer point at $0 \in \C$. Then,
	there is a $G_{\delta}$-dense set $\mathcal{U} \subset \diffalpha \times \diffalpha$ such that, whenever $(h_1, h_2)$
	lies in $\mathcal{U}$, the pseudogroup $\Gamma_{h_1,h_2}$ generated by $\tilde{f} = h_1^{-1} \circ f \circ h_1, \,
	\tilde{g} = h_2^{-1} \circ g \circ h_2$ on $D$ satisfies the following:
	\begin{enumerate}
		\item The stabilizer of every point $p \in D$ is either trivial or cyclic.
		
		\item There is a sequence of points $\{ Q_i\}$, $Q_i \neq 0$ for every $i \in \N^{\ast}$, converging to $0 \in \C$
		and such that every $Q_n$ is a hyperbolic fixed point of some element $W_i (\tilde{f}, \tilde{g}) \in \Gamma_{h_1,h_2}$. Furthermore the orbits under $\Gamma_{h_1,h_2}$ of $Q_{n_1}, \, Q_{n_2}$ are disjoint provided
		that $n_1 \neq n_2$.
	\end{enumerate}
\end{teo}

To conclude, we should only mention that Theorems~\ref{teo_diffeo1} and~\ref{teo_diffeo2} can be applied to much
larger classes of foliations. In fact, they can be applied to every class of foliations that are stable under perturbations
in the Krull topology such as, for example, those singularities whose resolution tree has only hyperbolic singular points.


\section{Integrability of foliations and related problems on \texorpdfstring{${\rm Diff}(\C^n,0)$}{dc0}}\label{sec_integrability}

In the context of singularities of holomorphic foliations in dimension~$2$, the topological nature of the foliation
and the existence of non-constant holomorphic first integrals possesses a surprisingly strong connection which was
put forward in the seminal paper \cite{MM}. In fact, the existence of first integrals for the foliations in question
can be read off as some clearly necessary topological conditions (recall that by ``first integral'' will always mean a non-constant first integral). More precisely, the following is proved in the mentioned
paper:

\begin{teo}[{\bf Mattei-Mossu Theorem}]\cite{MM}\label{teo_MM}
	Consider a holomorphic foliation $\fol$ defined on a neighborhood $U$ of the origin of $\C^2$. The foliation $\fol$ has a non-constant holomorphic first integral $f: U \to \C$ if and only if
	the following two conditions are satisfied:
	\begin{itemize}
		\item[1.] only a finite number of leaves accumulates on $(0,0)$;
		\item[2.] the leaves of $\fol$ are closed on $U \setminus \{(0,0)\}$.
	\end{itemize}
\end{teo}

It immediately follows from what precedes that the existence of a non-constant holomorphic first integral for a singular foliation
on $(\C^2,0)$ is a topological invariant. In other words:

\begin{corollary}\cite{MM}
	Consider two local foliations $\fol_1, \, \fol_2$ about $(0,0)
	\in \C^2$ that are topologically equivalent in the sense that there is a local homeomorphism $h$ around $(0,0) \in \C^2$ and taking
	the leaves of $\fol_1$ to the leaves of $\fol_2$. Then $\fol_1$ admits a non-constant holomorphic first integral if and only if so
	does $\fol_2$.	
\end{corollary}

\subsection{Extensions of Mattei-Moussu theorem to higher dimensions}

Possible generalizations of the above mentioned phenomenon have long attracted interest. First, a classical example
attributed to Suzuki and discussed by Cerveau and Mattei in~\cite{C-M} shows that the existence of meromorphic first integrals is no longer a
topological invariant. However, in dimension~$3$, many experts have wondered whether the existence of two ``independent"
holomorphic first integrals would constitute a topological invariant of the singularity. Recently, in \cite{P-R}, this
question was answered in the negative. Indeed, we proved the following:

\begin{teo}\cite{P-R}\label{contra-exemplo}
	Denote by $\fol$ and $\mcf$ the foliations associated to the vector fields $X$ and $Y$, respectively, given by
	\begin{align*}
		X &= 2xy \frac{\partial }{\partial x} + (x^3+2y^2) \frac{\partial }{\partial y} - 2yz \frac{\partial }{\partial z} \, , \\
		Y &= x(x - 2y^2 - y) \frac{\partial }{\partial x} + y(x - y^2 - y) \frac{\partial }{\partial y} - z(x - y^2 - y) \frac{\partial }{\partial z} \, .
	\end{align*}
	The foliations $\fol, \, \mcf$ are topological equivalent. Nonetheless $\fol$ admits two independent holomorphic
	first integrals while $\mcf$ does not.
\end{teo}

Our construction is inspired from Suzuki's example based on a simple observation that the existence of two independent
holomorphic first integrals may give rise to a meromorphic first integral for the restriction of the foliation to certain
invariant surfaces. In fact, let $\fol$ be a foliation on $\left(\C^3,0\right)$ admitting two (necessarily) non-constant
and independent holomorphic first integrals $F$ and $G$. Consider the decomposition of $F$ and $G$ into irreducible factors
\begin{eqnarray*}
	F & = & f_1^{m_1}\cdots f_k^{m_k} \\
	G & = & g_1^{n_1}\cdots g_l^{n_l} \, \, \, .
\end{eqnarray*}
Suppose that $F$ and $G$ have no common irreducible factor, modulo multiplication by nowhere vanishing functions. Then
the restriction of $G$ to, for example, $\{f_1=0\}$ is a non-constant holomorphic first integral for the restriction of
$\Fcal$ to the surface in question. In particular, the restriction of the foliation $\Fcal$ to $\{f_1 = 0\}$, viewed as a singular
foliation defined on a (possibly singular) surface, admits finitely many separatrices. In this case, all leaves of
$\Fcal|_{\{f_1 = 0\}}$ are ``fully identified" by $G$ in the sense that the restriction of $G$ to $\{f_1 = 0\}$ provides
a non-constant holomorphic first integral for $\Fcal|_{\{f_1 = 0\}}$. Assume now that $f_1$ is a common irreducible factor
for $F$ and $G$. Then the restrictions of both $F$ and $G$ to $\{f_1 = 0\}$ vanish identically. In this case, the leaves of
$\Fcal|_{\{f_1 = 0\}}$ cannot be distinguished by either $F$ or $G$. Nonetheless, it is possible to obtain a non-constant
first integral for the restriction of $\Fcal$ to $\{f_1 = 0\}$ as a function of $F$ and $G$. To be more precise, there exist
positive integers $n_1, m_1$ such that the function
\begin{equation}\label{eq1}
	\frac{F^{n_1}}{G^{m_1}}=\frac{f_2^{m_2n_1}\cdots f_k^{m_kn_1}}{g_2^{n_2m_1}\cdots g_l^{n_lm_1}}
\end{equation}
is a non-constant first integral of $\Fcal_{|_{\{f_1=0\}}}$. However, in general, this first integral is meromorphic rather
than holomorphic as shown by the simple example below.

\begin{example}
	Consider the holomorphic functions $F = xy$ and $G = xz$ which clearly define two independent holomorphic first integrals for
	the foliation associated to the vector field $X = x \partial/ \partial x - y \partial /\partial y - z \partial /\partial z$.
	Both $F, G$ vanish identically on the invariant manifold $\{x=0\}$. Nonetheless, the function $F/G = y/z$ provides a meromorphic
	first integral for the restriction of $\fol$ to this invariant manifold.
\end{example}

In view of the above observation, and recalling that the existence of a meromorphic first integral is not a topological
invariant, the definition of the foliations $\fol$ and $\mcf$ in Theorem~\ref{contra-exemplo} is itself inspired from the
Suzuki and Cerveau-Mattei examples in the following sense. The plane $\{z = 0\}$ is invariant by both $\fol, \, \mcf$ and
the restriction of $\fol$ (resp. $\mcf$) to this invariant manifold coincides with the foliation provided by Cerveau-Mattei
(resp. Suzuki). Furthermore $\fol$ and $\mcf$ were chosen so that the image of each leaf of $\fol$ (resp. $\mcf$) under the
projection map ${\rm pr}_2(x,y,z) = (x,y)$ is still a leaf of $\fol$ (resp. $\mcf$) and, in addition, a sort of ``saddle
behavior'' for their leaves with respect to the third axes was introduced (by ``saddle behavior" it is meant that as the
variable $x$ on the local coordinates of a leaf decreases to zero, the variable $z$ increases monotonically to exit a fixed
neighborhood of the origin). The ``saddle behaviour'' was carefully chosen so that the topological equivalence between the
restrictions of $\fol, \, \mcf$ to the invariant plane $\{z=0\}$ can be extended to an entire neighborhood of the origin.

From the above construction, the foliation $\fol$ possesses $\overline{F} = (y^2 - x^3)z^2$ and $\overline{G} = xz$
as independent holomorphic first integrals. Furthermore, it is clear that the foliation $\mcf$ cannot admit two independent
holomorphic first integrals. In fact, if this were the case, then the quotient between suitable powers of their first integrals
would provide a meromorphic first integral for the restriction of $\mcf$ to the invariant plane $\{z = 0\}$. As already mentioned,
it is known that such meromorphic first integral does not exist.

It can be noted that the singular set of the foliations considered in Theorem~\ref{contra-exemplo} is not reduced to
a single point and this might suggest that the ``correct'' generalization of Mattei-Moussu theorem involves isolated
singularities. This is actually not the case. As follows from the above construction, the existence of invariant surfaces over which the corresponding foliation is dicritical often constitutes an essential obstruction for the topological invariance of ``complete integrability''. Furthermore, in $3$-dimensional ambient spaces, there is vast evidence that completely integrable foliations with isolated singularities must admit an invariant surface over which the correspondent foliation is
dicritical. In fact, in the same paper, the following result was proved.

\begin{teo}\cite{P-R}\label{teodicritical}
	Let $\fol$ be a foliation by curves on $(\C^3,0)$ having an isolated singularity at the origin and admitting
	two independent holomorphic first integrals. Suppose that $\tilf$, the transform of $\fol$ by the one-point
	blow-up centered at the origin, has only isolated singularities which, in addition, are simple. Then $\fol$
	possesses an invariant surface over which the induced foliation is dicritical.
\end{teo}

Concerning the role played by the above mentioned invariant surfaces, recall that a deep study of topological properties
of foliations on $(\C^2,0)$ possessing meromorphic first integrals was conducted by M. Klughertz in~\cite{Martine}. Her
techniques yield several examples where ``topological invariance" for the existence of meromorphic first integrals fails.
Relatively simple adaptations of the proof of Theorem~\ref{contra-exemplo} then enable us to obtain several other examples
of foliations on $(\C^3,0)$ for which the ``topological invariance" of the existence of two independent holomorphic first
integrals is not verified. We conjecture, however, that if we are given two foliations by curves on $(\C^3,0)$, $\fol_1, \,
\fol_2$, that are topologically equivalent and do not admit invariant surfaces over which the induced foliations are dicritical,
then $\fol_1$ admits two holomorphic first integrals if and only if so does $\fol_2$.

Recall that the Seidenberg desingularization theorem plays a major role in the study of singular foliations in dimension~$2$
and, in particular, it is used in the topological characterization of integrable foliations. A completely faithful generalization
of the Seidenberg result for foliations on $3$-manifolds cannot exist, since some non-simple singularities are persistent under
blow-ups (cf. Section~\ref{Sec:resolution}). Nonetheless final models on a desingularization process of foliations on $3$-manifolds
have been described on different papers such as \cite{C-R-S}, \cite{MQ-P}, \cite{P} and \cite{RR_Resolution} (see Section~\ref{Sec:resolution} for details). In a first moment,
the idea to remove the generic condition from Theorem~\ref{teodicritical}, and/or to eventually prove the above mentioned conjecture,
consists of showing that this ``special type'' of singular points cannot appear in the desingularization procedure of $\fol$ provided
that $\fol$ is completely integrable. Following some discussions with D. Panazzolo, this assertion can probably be established by
building on the material of the mentioned above papers.

\bigbreak

Another celebrated theorem by Mattei and Mossu in~\cite{MM} can be stated as follows

\begin{teo}\cite{MM}\label{teo:formali_first_int}
	Let $\fol$ be a codimension-$1$ holomorphic foliation on $(\C^n, 0)$. If $\fol$ admits a formal first integral, then $\fol$ admits a holomorphic first integral.
\end{teo}

Their work on the existence of first integrals has also motivated Malgrange's Theorem in~\cite{malgrange}. The relationship between formal and holomorphic first integrals for higher codimension foliations remains, however, quite mysterious so that it is natural to begin with the case of $1$-dimensional foliations on $(\C^3,0)$. This is the simplest case outside the reach of Mattei-Moussu's results. In this context, Cerveau asked whether a $1$-dimensional (holomorphic) foliation on $(\C^3,0)$ admitting one (resp. two independent) formal first integral(s) must possess holomorphic one (resp. two independent) first integral(s) as well. This question has been answered in the negative for the case of one formal fist integral. In fact, in a joint work with A. Belotto, M. Klimes and J. Rebelo, the following was proved:

\begin{teo}\cite{BKRR}\label{teo:conterexample_Cerveau}
	Consider the family $X_{a,b,c}$ of vector fields on $\C^3$ defined by
	\begin{equation}\label{eq:vf}
		X_{a,b,c}= x^2 \frac{\partial}{\partial x} + (1+ax)\left[ y_1\frac{\partial}{\partial y_1}-
		y_2 \frac{\partial}{\partial y_2} \right] +  bxy_2\frac{\partial}{\partial y_1} + cxy_1\frac{\partial}{\partial y_2} \, ,
	\end{equation}
	where $a$, $b$, and $c$ are complex parameters. Assume that the parameters  are such that
	\[
	\cos (2\pi a) \neq \cos (2\pi\sqrt{a^2+bc}) \, .
	\]
	Then the vector field $X_{a,b,c}$ possesses no (non-constant) holomorphic first integrals, albeit it does
	possess formal first integrals.
\end{teo}

In particular, the vector field $X_{1,1,1}$ obtained by setting $a = b = c = 1$ admits a formal first integral but no holomorphic one.

The proof of Theorem~\ref{teo:conterexample_Cerveau} relies on the standard theory of linear systems (normal forms and Stokes phenomena among others). The example in the mentioned theorem was constructed after having computed (through the same techniques as those used in the proof of this theorem) the normal form and the corresponding Stokes phenomena of the vector field
\[
Y_A= - \frac{1}{2} x^4 \frac{\partial}{\partial x} + \left(z-\frac{1}{2}x^3y\right)
\frac{\partial}{\partial y} + (y- x^3z) \frac{\partial}{\partial z} \,,
\]
which corresponds to a saddle-node singularity appearing in a convenient birational model for the compactified Airy equation. The Airy equation is usually seen as a toy model for Painlev\'e II equation. We have checked that $Y_A$ admits a formal ``meromorphic'' first integral (i.e., there is a formal first integral taking on the form $F/G$, where $F, \, G \in \C[[x, y, z]]$), but it admits no holomorphic or meromorphic first integral.

The question if a holomorphic vector field admits two independent formal first integrals (i.e. two formal first integrals $F, \, G$ such that $dF \wedge dG \not \equiv 0$) also admits two independent holomorphic first integrals remains unknown.


\subsection{Problems on \texorpdfstring{${\rm Diff} \, (\C^n,0)$}{dc0} related with first integrals for foliations}

Motivated by the above examples, we conducted a more in-depth study of Mattei-Moussu's results as well as their possible
generalizations. There, it should be noted that the Mattei-Moussu argument also states that the existence of these first
integrals can be detected at the level of the topological dynamics associated with the holonomy pseudogroup of the foliation
in question. Let us make it precise.

Consider a holomorphic foliation $\fol$ on $(\C^2,0)$ and assume it admits a holomorphic first integral. Then so does the
holonomy pseudogroup of the foliation in question. More precisely, the holonomy pseudogroup of the given foliation corresponds
(up to a change of coordinates) to a group of rotations being, in particular, finite. It then follows that every element of the
mentioned pseudogroup has finite orbits. In fact, the pseudogroup itself has finite orbits. Recall that

\begin{defi}
	We say that a group $G$ has finite orbits, if there exists a sufficiently small neighborhood $V$ of the origin such that the set $\mathcal{O}_V^G(p)$ is
	finite for every $p \in V$, where
	\[
	\calO_V^G (p) = \{q \in V \; \, : \; \,  q = h(p), \; h \in G \; \; {\rm and} \; \; p \in {\rm Dom}_V (h) \} \, .
	\]
\end{defi}

The central point of the proof of Mattei-Moussu theorem is a converse for the previous statement which is valid for subgroups of
$\diff$, namely:

\begin{prop}\label{prop_finitude}\cite{MM}
	Let $G$ be a finitely generated pseudogroup of $\diff$. Assume that $G$ has finite orbits. Then $G$ is itself finite.
\end{prop}

The proof of the Mattei-Moussu's Theorem (i.e. Theorem~\ref{teo_MM}) goes essentially as follows. From the proposition above, it follows that the holonomy pseudogroup of $\fol$ is conjugate to a finite group of rotations (being finite it is actually cyclic). It immediately follows that it admits a first integral of the form $z \mapsto z^n$, for a certain $n \in \N^{\ast}$. Then, they proceed to extends this first integral along the leaves of the foliation to derive a holomorphic first integral for the foliation as well.

With respect to the extension of the first integral of the holonomy pseudogroup through the saturated of leaves, the following should
be noted. The resolution of a foliation as $\fol$ is such that every singular point in the final model is in the Siegel domain (in other words,
if $\lambda_1, \lambda_2$ are the eigenvalues at a singular point then both eigenvalues are different from zero and $\lambda_1/\lambda_2
\in \R^-$, cf. Section~\ref{sec:basics}). One such singular point is such that the corresponding foliation possesses exactly two separatrices. Mattei proved, in an
unpublished manuscript, the following:

\begin{prop}[\cite{Mattei_unpublished}]\label{prop:Mattei_unpub}
	Let $\fol'$ be a foliation on $(\C^2,0)$ with a singular point of Siegel type at the origin. The saturated of a local transverse section to any one of the two separatrices,  together with the other separatrix contains a neighbourhood of the singularity in question.
\end{prop}

In the above context, Proposition~\ref{prop:Mattei_unpub} plays a fundamental role related to showing that the closure of the saturated by $\fol$ of the domain of the initial first integral $z \mapsto z^n$ actually contains a neighborhood of the mentioned singular point.

\bigbreak

Extensions of all of the previous results to foliations on higher dimensional manifolds were provided, at least under suitable
conditions. Let us start by stating the extention obtained with respect to this last result, that is the result ensuring that
the saturated of a transversal section to a separatrix through a singular point in the Siegel domain, together with the other
separatrix, constitutes a neighborhood of the origin. In higher dimensions, the result becomes:

\begin{prop}\cite{Reis06,RR_LN}\label{prop_gener_Mattei}
	Let $\fol$ be a singular foliation associated to a holomorphic vector field $X$ with an isolated singularity at the origin of $\C^n$.
	Suppose that the origin belongs to the Siegel domain and satisfy the following conditions:
	\begin{itemize}
		\item[(a)] The eigenvalues $\dl_1, \ldots, \dl_n$ of the linear part of $X$ at $0 \in \C^n$ are all different from zero and there exists
		a straight line through the origin, in the complex plane, separating (for example) $\dl_1$ from the remainder eigenvalues.
		
		\item[(b)] Up to a change of coordinates, $X = \sum_{i=1}^n \dl_ix_i(1+f_i(x)) \partial /\partial x_i$, where $x=(x_1,\ldots,x_n)$ and
		$f_i(0)=0$ for all $i$.
	\end{itemize}
	Then the saturated of a transversal section to the separatrix associated with the eigenvalue $\dl_1$ at a point sufficiently close to the
	origin, together with the invariant manifold transverse to the mentioned separatrix contains a neighborhood of the origin.
\end{prop}

To begin with, it should be noted that if $X$ is a vector field on $\C^3$ with an isolated singularity at the origin and of ``strict Siegel
type'' (i.e. the convex hull of $(\dl_1,\dl_2, \dl_3)$ contains a neighbourhood of the origin), then conditions (a) and (b) are immediately
satisfied (cf. \cite{CKP} for item (b); with respect to item (a), it is clear there exists at least one eigenvalue $\dl_i$ such that the
angle between $\dl_i$ and the other eigenvalues is greater than $\pi/2$). This result is then ``generic'' in dimension~$3$. There are,
however, examples of vector fields whose origin belongs to the boundary of the convex hull of $(\dl_1, \dl_2, \dl_3)$ does not satisfying
item~(b) (cf.~\cite{canille}). Nonetheless, the existence of the three invariant hyperplanes plays a role in the proof of the theorem in question.

Let us say a few words about the proofs. In dimension~$2$ the proof is based on the following. Consider standard coordinates $(x_1,x_2)$
where $X$ takes on the form of item (b). Fix the separatrix $S$ given in the present coordinates as the $x_1$-axis and let $\Sigma$ stand
for a transversal section to $S$ at a point sufficiently close to the origin. Assume, without loss of generality, that $\Sigma \subseteq
\{x_1 = \varepsilon\}$ for some arbitrarily small $\varepsilon \in \R$. Then
\begin{itemize}
	\item consider on $S$ the loop given as $x_{1,l}(t) = \varepsilon e^{2\pi i t}$, with $0 \leq t \leq 1$. A ``kind'' of solid torus around
	the origin can be obtained by taking the lift of the loop $x_{1,l}$ for all points on $\Sigma$;
	
	\item next, consider on $S$ the ``radial directions'' given by $x_{1,A}(t) = Ae^{-t}$, $t>0$, for every $A \in \C$ with $|A| = \varepsilon$,
	and take the lift of the mentioned path through every single point of the previous ``solid torus''.
\end{itemize}
So, fixed a point $p$ on the above ``solid torus'', let $A = pr_1(p)$, where $pr_1$ stands for the projection map of $\C^2$ on the first
component, i.e. in local coordinates $(x,y)$ we have $p_r(x,y) = x$. By providing precise estimates, Mattei proved that if $\varphi(t) =
(x_{1,A}(t),y(t))$ stands for the lift of the path $(x_{1,A}(t),0)$ through $(A,0)$ along the leaf passing through $p$, then $y(t)$ escapes
from any small neighborhood of the origin. Essentially, the lift of the mentioned path has a saddle behavior and it is this saddle behavior
that ensures that the saturated of the transversal section must contain a neighborhood of the origin, up to joining the other separatrix.

In higher dimensions, extra conditions had to be imposed to ensure that the leaves of the foliation present a similar behavior. Before accurately indicating
the role of the imposed conditions in the proof of Proposition~\ref{prop_gener_Mattei}, let us make some remarks based on a concrete example.

\begin{example}
	{\rm Consider a vector field $X$ on $(\C^3,0)$ taking on the form $X = \sum_{i=1}^3 \dl_ix_i(1+f_i(x)) \, \partial /\partial x_i$ with
		$(\dl_1,\dl_2,\dl_3) = (1,1+i,-2-i)$, so that the origin is a singular point of strict Siegel type. Furthermore there exists a straight
		line through the origin, in the complex plane, separating $\dl_1$ from the remainder eigenvalues. Therefore, Proposition~\ref{prop_gener_Mattei}
		states that the saturated of a transversal section to the separatrix $S$ given as the $x_1$-axis through a point arbitrarily close to the
		origin, jointly with the invariant hyperplane $\{x_1=0\}$, constitutes a neighborhood of the origin. It can easily be checked however that
		the lift of the ``radial directions'' given by $x_{1,A}(t) = Ae^{-t}$, $t>0$, for every $A \in \C$ with $|A| = \varepsilon$, through any
		point along $\Sigma_A = \{x_1 = A\}$ does not present a ``saddle behaviour''. To be more precise, if $\varphi(t) = (x_{1,A}(t),x_2(t),
		x_3(t))$ stands for the lift of the path $(x_{1,A}(t),0,0)$ through $(A,0,0)$ along the leaf passing through $p \in \Sigma_A$, then
		$x_2(t)$ goes to zero as $t$ goes to infinity.}
\end{example}

The fact that $|x_2(t)|$ goes to zero as $t$ goes to infinity, in the previous example, happens since the angle between $\dl_1$ and $\dl_2$ is strictly less than $\pi/2$.
In any case, we were able to establish the saddle behavior by taking the lift along paths distinct from the ``radial lines'' on $S$ but
still accumulating at the origin. Let us precise the path in question.

\begin{figure}[hbtp]
	\centering
	\includegraphics[scale=0.65]{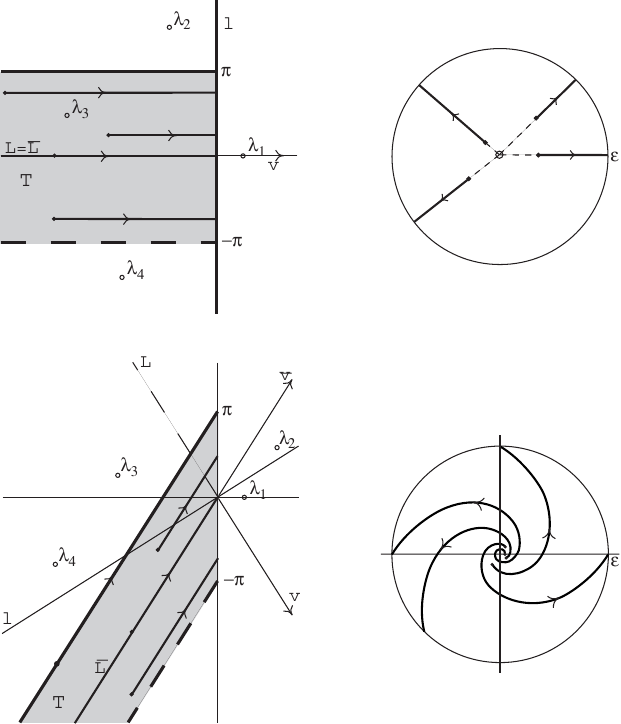}
	\caption{}\label{fig1}
\end{figure}

Assume without of generality that $\dl_1 = 1$. So, in the complex plane, let $l$ be a straight line through the origin separating
$\dl_1$ from the remaining eigenvalues. Consider then the straight line orthogonal to $l$ at the origin and denote be $L$ the part
of this straight line that is contained in the left half-plane (negative real part). Finally, denote by $\bar{L}$ the complex conjugate
of $L$. Suppose that $v = \da + i \db$, with $\da > 0$, is a directional vector of $L$. Then let
\[
T = \{z \in \C \, : \, z = x + iy , x \in \bar{L}, -\pi < y \leq \pi\}
\]
(cf. figure \ref{fig1}). It can easily be checked that the image of T by the application map $\phi: \, z \to \de e^{z}$ covers
$\{z: |z| \leq \de\} \setminus \{0\}$, being the map in question one-to-one. Moreover, the image by $\phi$ of the elements
$z = iy$, with $-\pi < y \leq \pi$, corresponds to the circle in $S$ of radius $\de$ and centered at the origin.

For every $y \in ]-\pi, \pi]$ we fix, let $c_y(t)$ be the path on the complex plane defined by
\[
c_y(t) = iy + \frac{1}{v} \, t \, ,
\]
with $v$ as above, for $t \in ]-\infty,0]$. Consider the logarithmic spiral curve contained in the $x_1$-axis given by $r_y (t) =
(\de e^{c_y(t)}, 0 , \ldots, 0)$ for $t \in ]-\infty,0]$. The spiral curve is such that $|\de e^{c_y(0)}| = \de$ and $\de e^{c_y(t)}$
goes to zero as $t$ goes to $-\infty$. Fix now and element $z \in \{x_1 = \de e^{iy}\}$ and let $r_z$ be the lift of $r_y$ through
the lift through $z$. The corresponding lift has a saddle behavior in the sense that the modulus of every component of $r_y$ increases
as the modulus of the first component decreases. In fact, if $\varphi(t) = (\de e^{c_y(t)}, x_2(t), \ldots, x_n(t))$ stands for the
mentioned lift, then $x_2(t), \ldots, x_n(t)$ satisfies
\[
\begin{cases}
	\frac{dx_2}{dt} = \frac{\dl_2}{v} x_2 \left( 1 + A_2 (\de e^{c_y(t)},x_2,\ldots,x_n) \right)\\
	\quad \vdots\\
	\frac{dx_n}{dt} = \frac{\dl_n}{v} x_n \left( 1 + A_n (\de e^{c_y(t)},x_2,\ldots,x_n) \right)
\end{cases}
\]
for some holomorphic functions $A_2, \ldots, A_n$. It is the fact that the angle between $v$ and $\dl_i$ is greater than $\pi/2$ that
ensures that $|x_i|$ increases as $t$ goes to $-\infty$ (cf. \cite{RR_LN} for precise calculations on the estimates).

More recently, in his thesis \cite{Chaves_thesis} (see also~\cite{Chaves_paper}), F. Chaves unifies and significantly extends the previous results. His theorem is essentially sharp and I will quickly review his theorem in the sequel so as to be able to comment on the progress made compared to previous works.

First, Chaves works with what he calls crossing type foliations, these are foliations possessing a
(smooth) invariant curve as well as an invariant hyperplane transverse to the mentioned invariant
curve. This assumption is very natural in view of basic issues of the problem and, actually, it is
weaker than the analogous versions used by Elizarov-Il'yashenko \cite{E-IY} and by myself. He also assumes
that the eigenvalue associated with the direction of the invariant curve is non-zero which is, indeed, a
much weaker condition than the Siegel type singularity condition used in the mentioned previous
works (more on this below). Finally, he assumes a very weak non-resonance condition involving
the eigenvalues associated to directions tangent to the invariant hyperplane; weak as it is, this
condition turns out to be necessary as shown by Chaves in his thesis work.

Under the preceding assumptions, he shows that two foliations $\fol$ and $\mathcal{G}$ on $(\C^n, 0)$ having conjugate
linear parts and conjugate holonomy maps arising from the invariant curves are necessarily
analytically equivalent.

Proposition~\ref{prop_gener_Mattei} becomes a Corollary of Chaves' result in the sense that his statement weakens every single condition imposed in Proposition~\ref{prop_gener_Mattei}. Moreover, the assumption about isolated singular points is not needed anymore. The main advantages of Chaves theorem can be summarized as follows:
\begin{itemize}
	\item In all previous works, the foliation was always assumed to have only non-zero eigenvalues (lying in the Siegel
	domain). In particular, the singular point was isolated. This condition is dramatically
	dropped by Chaves that only requires a very minor non-resonance condition. In this way, his
	results applies to saddle-node singularities including non-isolated ones. Already in dimension
	2, his theorem provides a very welcome alternative proof of a classical theorem due to
	Martinet and Ramis.
	
	\item Also even in the case where all eigenvalues are different from zero (isolated singularity in
	the Siegel domain), Chaves's theorem entails a major progress for $n \geq 5$. This is down
	to the assumption on the ``distinguished eigenvalue'' that can be separated from the other
	eigenvalues by a straight line. To understand the issue, it is enough to think of 5 eigenvalues
	defining a (regular) pentagon in $\C^{\ast}$. In this case a ``distinguished eigenvalue'' does not
	exist and the situation is stable under perturbation of the eigenvalues. Thus, even in the Siegel
	case, there are open sets of possibilities that are not covered by the former theory while they
	are covered by Chaves' theorem, at least up to the very minor non-resonance condition.
	
	\item Finally, whereas the non-resonance condition imposed by Chaves is a very minor one, he
	can still prove that it is necessary. In fact, when this condition is dropped, he provides
	several examples of non-analytically equivalent pairs of foliations satisfying all the remaining
	conditions. This is the ultimate confirmation of the sharpness and depth of his methods.
\end{itemize}

I might still add a general remark to the above discussion. The method used by Chaves does not
rely on the ``path lifting technique'' used by Mattei-Moussu, Elizarov-Il'yashenko, and by myself.
Instead, he builds on an idea of A. Diaw in \cite{Arame_thesis} and \cite{Arame_paper}. Basically this technique consists of first defining the
equivalence between the foliations over the product of an annulus with a small disc (thus avoiding
any difficulty arising from the singular point) and, in a second moment, of extending only the
conjugating map to a neighborhood of the origin.

\bigbreak

Let us go back to Proposition~\ref{prop_finitude} and extensions of it. This proposition is, in fact, quite important and similar characterizations of pseudogroups having finite orbits
(or more generally locally finite orbits) has many applications. To begin with, it should be noted that the statement of the proposition in question is no longer valid in higher dimensions. The simplest example of an element of ${\rm Diff}\, (\C^2,0)$ having finite orbits is obtained by setting $F(x,y) = (x + f(y),y)$ with $f(0) = 0$. For example, by letting $f(y) = 2 \pi iy$, the resulting diffeomorphism $F$ can be realized as the holonomy map of the foliation associated with the vector field
\[
X = y\frac{\partial}{\partial x} + z\frac{\partial}{\partial z} \, ,
\]
with respect to the invariant curve $\{x = y = 0\}$. In this example, however, the reason for the element $F$ to have orbits relates to the fact that the non-fixed points of $F$ leave every single neighborhood of the origin after a finite number of iterations and not from the fact that that $F$ is periodic. In \cite{RebeloR_note}, the finiteness of a
pseudogroup having finite orbits was established in $\diffn$ under rather restrictive conditions involving isolated fixed points.
Also several non-trivial examples of infinite pseudogroups having {\it finite}\, orbits were provided already
in ${\rm Diff}\, (\C^2,0)$. These examples are more interesting than the previous one and, many of them, associated with singularities in the Siegel domain.

Yet, a far more relevant question was to figure out, in the general case,
which type of algebraic conditions a subgroup of $\diffn$ possessing finite orbits should verify. Differential
Galois theory, as well as Morales-Ramis theory cf. \cite{MRS}, suggests that a pseudogroup having finite orbits may be (virtually)
solvable. In the paper~\cite{RR}, we confirmed this suggestion for subgroups of $\diffd$. More precisely,
the following was proved:

\begin{teo}\label{teo_solvable}\cite{RR}
	Suppose that $G$ is a finitely generated pseudosubgroup of $\diffd$ with locally discrete orbits. Then $G$ is virtually solvable.
\end{teo}

Let us recall what we mean by a virtually solvable group.

\begin{defi}
	A group $G$ is said to be virtually solvable if it contains a normal and solvable subgroup $G_0$ of $G$ of finite index.
\end{defi}

\begin{rem}
	{\rm It should be noted that the condition ``virtually'', on the statement of Theorem~\ref{teo_solvable}, is natural in the sense that the assumption of having locally discrete orbits is stable under finite extensions of the group. Furthermore, this condition is also necessary. In fact, consider $A_5$, the group of even permutation on ~$5$ elements, that that can be realized as a subgroup of ${\rm PSL} \, (2, \C)$ and, consequently, of $\diffd$. This group is finite and, consequently, it is a group with locally discrete orbits. However, $A_5$ is not solvable.}
\end{rem}

Since, from the point of view of differential Galois theory, solvable groups are associated with systems ``integrable by quadratures'', Theorem~\ref{teo_solvable} states that the
``integrability of the group'' can be detected through its topological dynamics. Note that the assumption of having
locally discrete orbits (as opposed to finite orbits) was introduced so as to allow more general types of first integrals
including meromorphic ones.

Whereas the statement of Theorem~\ref{teo_solvable} appears to be sharp, it can considerably be strengthened in the case of groups of diffeomorphisms tangent to the identity. To state a sharper result on these groups, let us recall the notion of {\it recurrent point}. Let then $\diffdd$ denote the subgroup of $\diffd$ consisting
of diffeomorphisms tangent to the identity.

\begin{defi}
	Let $G$ be a subgroup of $\diffdd$. A point $p$ is said to be {\it recurrent}\, if its orbit
	under $G$ accumulates non-trivially on $p$ itself .
\end{defi}

Note that the definition excludes points having ``periodic'' orbit. The following result provides then strong quantitative information concerning the set of recurrent points.

\begin{teo}\label{teo_non_solvable}\cite{RR}
	Suppose that $G \subseteq \diffdd$ is non-solvable. Then there exists a neighborhood of the origin $U$
	and a countable union $K \subset U$ of proper analytic subsets of $U$ such that every point in $U \setminus K$ is recurrent
	for $G$ (in particular the set of recurrent points has full Lebesgue measure).
\end{teo}

Theorem~\ref{teo_non_solvable} claims, in particular, that a non-solvable subgroup of $\diffdd$ has plenty of points whose orbit is not locally discrete. Moreover, we manage to show that for a ``generic group'' (in a sense precised in the paper) the set of non-recurrent points is reduced to the origin $\{(0,0)\}$.

The results obtained in this paper had later been generalized by Ribon to higher dimensions (see~\cite{Ribon}).

\vspace{0.25cm}

\section{Resolution of singularities for 1-dimensional foliations and for vector fields on \texorpdfstring{$3$}{3}-manifolds}\label{Sec:resolution}

In this section it will be presented the main results obtained in the paper~\cite{RR_Resolution} along with the basic notions
needed to make their statements intelligible. A discussion about the place of these results in the current state-of-art in the
area will also be presented.

Recall first that a singular, one-dimensional holomorphic foliation $\fol$ on $(\C^n,0)$ is nothing but the (singular) foliation
defined by the local orbits of a holomorphic vector field defined on a neighborhood of the origin and having zero-set of codimension
at least~$2$. A simple consequence of Hilbert nullstellensatz is that, up to multiplying vector fields by a meromorphic function,
every meromorphic vector field $X$ on $(\C^3,0)$ induces a singular holomorphic foliation on a neighborhood of the origin. This
foliation will be called the foliation associated with $X$. Clearly two (meromorphic) vector fields have the same associated foliation
if and only if they differ by a multiplicative (meromorphic) function. Conversely, a vector field $X$ inducing a given foliation $\fol$
will be called a representative of $\fol$ if $X$ is holomorphic and the set of zeros of $X$ has codimension at least two. In other words,
a representative vector field of $\fol$ is any holomorphic vector field tangent to $\fol$ and having a zero-set of codimension
at least~$2$.

There follows from the preceding that there is no point in considering ``singular meromorphic foliations'' since all foliations in this
category would, in fact, be holomorphic. Similarly, (singular) holomorphic foliations have empty zero-divisor since their singular sets
have codimension at least~$2$. In other words, whenever we are exclusively concerned with foliations, we can freely eliminate any
(meromorphic) common factor between the components of a vector field tangent to the foliation to obtain a representative vector field.
Naturally this cannot be done if we are focusing on an actual fixed vector field $X$ as it so often happens.

In the above mentioned context of singular points, {\it resolution theorems} - also known as {\it desingularization theorems} -
are geared towards foliations in that we are ``free'' to eliminate non-trivial common factors between the components of a vector
field whenever these common factors arise from transforming a representative vector field by a birational map. To further clarify
these issues, we may recall that the prototype of all ``resolution theorems'' for foliations is provided by Seidenberg's theorem,
which is valid for foliations defined on a two-dimensional ambient space, namely we have:

\begin{teo}[{\bf Seidenberg Theorem}]\cite{seiden}
	Let $\fol$ be a singular
	holomorphic foliation defined on a neighborhood of $(0,0) \in \C^2$. Then, there exists a finite sequence of blow-up maps, along with transformed foliations $\fol_i$ ($i=1, \ldots, n$)
	\[
	\fol = \fol_0 \stackrel{\Pi_1}\longleftarrow \fol_1 \stackrel{\Pi_2}\longleftarrow \cdots
	\stackrel{\Pi_l}\longleftarrow \fol_n
	\]
	such that the following holds:
	\begin{itemize}
		\item Each blow-up map $\Pi_i$ ($i=1, \ldots, n$) is centered at a singular point of $\fol_{i-1}$.
		
		\item All singular points of $\fol_n$ are {\it elementary}, i.e. $\fol_n$ is locally given by a representative
		vector field $X_n$ whose linear part at the singular point in question has at least one eigenvalue different from zero (cf. below).
	\end{itemize}
\end{teo}

Whereas Seidenberg's theorem is directly concerned with foliations, it is also very effective when applied to vector fields defined
on complex surfaces.
The general principle to use Seidenberg theorem to study vector fields - as opposed to foliations - consists of
applying Seidenberg theorem to the associated foliation while also keeping track of the divisor of zeros/poles of the
transformed vector field. In line with this point of view, Seidenberg's theorem is equally satisfying: the structure of the resolution map
(the composition of the blow-ups $\Pi_i$) is such that the transform of holomorphic vector fields retains its holomorphic
character (here the reader is reminded that the transform of a holomorphic vector field by a birational map is, in general,
a meromorphic vector field). More generally, Seidenberg's procedure allows for an immediate computation of
the zero-divisor of the transformed vector field. For example, if we blow-up a vector field $X$ having an isolated
singularity at $(0,0) \in \C^2$ and denote by $k$ the degree of the first non-zero homogeneous
component of the Taylor series of $X$ at $(0,0)$, then the zero-divisor of the blow-up of $X$ coincides with
the exceptional divisor and has multiplicity $k-1$ (unless the first non-zero homogeneous component of $X$ is actually a multiple of the radial vector field
- $R = x\partial /\partial x + y \partial /\partial y$ - in which case the multiplicity is~$k$).

In dimension~$2$, the classical Seidenberg Theorem provides an optimal algorithm for simplifying the singularities of a foliation.
On the other hand, when one moves to dimension~$3$, the situation is no longer so simple. The well-known example of Sancho and Sanz
shows the existence of foliations in $(\C^3,0)$ that cannot be reduced by standard blow-up centered at the singular set of the
foliation in question. To be more precise, they have shown that the foliation associated with the vector field
\[
X = x\left( x \frac{\partial}{\partial x} - \alpha y \frac{\partial}{\partial y}  - \beta z \frac{\partial}{\partial z}\right)
+ xz \frac{\partial}{\partial y} + (y-\lambda x) \frac{\partial}{\partial z}
\]
possesses a strictly formal separatrix $S = S_0$ through the origin such that the singular point $p_n$ (selected by the transformed
separatrices $S_n$ in the sense that they correspond to the intersection of $S_n$ with the excetional divisor) is a nilpotent singular
point for the corresponding foliation, for all $n \in \N$. In fact, the representative of the singular point $p_n$ is given by a vector
field on the above $3$-parameter family. The fact that every separatrix $S_n$ is stricty formal says that even in the case we allow
blow-ups to be centered at analytic invariant curves that {\it are not}\, necessarily contained in singular set of the foliation, a resolution
procedure still does not exist.

The generalization of Seidenberg's theorem to foliations on $(\C^3,0)$ is a very subtle problem. A very satisfactory answer is
provided in \cite{MQ-P-preprint}, \cite{MQ-P} and it relies heavily on a previous result by Panazzolo in \cite{P}. Slightly later,
the topic was revisited from the point of view of valuations in \cite{C-R-S}. The ``final models'' in the resolution theorem proved
in \cite{C-R-S} are, however, not as accurate as those in \cite{MQ-P}. The paper \cite{RR_Resolution} grew out of an attempt to use
the mentioned results to obtain a sharper resolution result which {\it would hold for the special class of holomorphic foliations}\,
which is associated with semicomplete vector fields (cf. Theorem~\ref{teo:B} below). Whereas the class of foliations associated with semicomplete
vector fields is rather special, it contains the underlying foliations of all complete vector fields as well as many foliations arising
in the context of Mathematical Physics and the importance of these examples justifies the interest in a sharper (or ``simpler'')
resolution statement valid only for this class of foliations.

The resolution theorem in~\cite{C-R-S} was not really suited to our needs because the corresponding ``final models'' were not accurate
enough. As to the resolution theorem in \cite{MQ-P}, we were unsure of the behavior of vector fields - as opposed to foliations - under
their procedure. Basically, we did not know if the weighted blow-ups on Panazzolo's algorithm \cite{P} always transform holomorphic vector fields
on holomorphic vector field rather than meromorphic ones (question that does not arise in the context of foliations, as explained above).
In fact, it is convenient to point out that, in full generality, the transform of a holomorphic vector field by a birational map is a meromorphic vector
field. To provide an explicit example.

\begin{example}
	{\rm Consider the holomorphic vector field $X = F(x,y,z) \partial /\partial x + G (x,y,z) \partial /\partial y + H (x,y,z)
		\partial /\partial z$ where $F(x,y,z) =y$ and $G$ and $H$ are such that the $z$-axis $\{ x=y=0\}$ is contained in the singular set of
		$X$. Let $(x,t,z)$ be coordinates for the weighted blow-up (of weight~$2$) centered at the $z$-axis in which the corresponding projection
		map $\Pi$ is given by $\Pi (x,t,z) = (x^2, tx, z)$. A direct inspection shows that the corresponding transform $\Pi^{\ast} X$ of $X$ is
		given by
		\begin{eqnarray*}
			\Pi^{\ast} X & = & \frac{1}{2x} F(x^2, tx, z) \frac{\partial}{\partial x} + \left[ -\frac{t}{2x^2} F(x^2, tx, z) +
			\frac{1}{x} G(x^2, tx, z) \right]\frac{\partial}{\partial t} + \\
			& & \, + H (x^2, tx, z) \frac{\partial}{\partial z} \, .
		\end{eqnarray*}
		Clearly $F(x^2, tx, z)/2x$ and $G(x^2, tx, z) /x$ are both holomorphic but $t F(x^2, tx, z) /2x^2$ is {\it strictly meromorphic}.
		Therefore $\Pi^{\ast} X$ is meromorphic with poles over the exceptional divisor.}
\end{example}

Although checking whether or not Panazzolo's algorithm in \cite{P} is such that the transforms of holomorphic vector fields retain their
holomorphic character should be straightforward, the algorithm itself is rather involved with many different cases so that we were very grateful
to the referee of our paper~\cite{RR_Resolution} for confirming that this is, in fact, the case. In other words, holomorphic vector fields
are transformed into holomorphic vector fields by the algorithm in \cite{P}. Still,  when studying the papers in question, we felt it would be nice to try and
complete the work of Cano-Roche-Spivakovsky~\cite{C-R-S} by deriving ``final models'' similar to those of \cite{MQ-P},
which are described in Theorem~\ref{teo:A} below.

\begin{teo}\cite{RR_Resolution}\label{teo:A}
	Let $\fol$ denote a (one-dimensional) singular holomorphic foliation defined on a neighborhood of $(0,0,0) \in \C^3$. Then there
	exists a finite sequence of blow-up maps along with transformed foliations
	\begin{equation}
		\fol = \fol_0 \stackrel{\Pi_1} \longleftarrow \fol_1 \stackrel{\Pi_2} \longleftarrow \cdots
		\stackrel{\Pi_l}\longleftarrow \fol_n \label{blowup-resolution_globaldesingularization}
	\end{equation}
	satisfying all of the following conditions:
	\begin{itemize}
		\item[(1)] The center of the blow-up map $\Pi_i$ is (smooth and) contained in the singular set of $\fol_{i-1}$, $i=1, \ldots , n$.
		
		\item[(2)] The singularities of $\fol_n$ are either elementary or persistent nilpotent singular points.
		
		\item[(3)] The number of persistent nilpotent singularities of $\fol_n$ is finite and each of them can be turned into
		elementary singular points by performing a single weighted blow-up of weight~$2$.
	\end{itemize}
\end{teo}

Recall that a nilpotent singular point $p$ is a point where the foliation admits a representative holomorphic vector field $X$ such that $X(p) = 0$ and $DX(p)$ is nilpotent but non-zero. A persistent nilpotent singular point is, roughly speaking, a nilpotent singular point for which any sequence of admissible standard blow-ups (i.e. standard blow-ups whose center is contained in the singular set of the foliation) always possess a nilpotent singular point. The precise definition appears in Section~$4$ of~\cite{RR_Resolution} and a normal form for these singularities is provided by Proposition~$3$ of the same paper. Let us recall their normal form.

\begin{prop}\cite{RR_Resolution}
	Let $\fol$ be a singular holomorphic foliation defined on a neighborhood of the origin of $\C^3$ and assume that the origin is a persistent nilpotent singularity of $\fol$. Then, up to finitely many one-point blow-ups, there exist local coordinates and a holomorphic vector field $X$ representing $\fol$ and having the form
	\[
	(y + f(x,y,z)) \frac{\partial}{\partial x} + g(x,y,z) \frac{\partial}{\partial y} +
	z^n \frac{\partial}{\partial z}
	\]
	for some $n \geq 2 \in \N$ and some holomorphic functions $f$ and $g$ of order at least~$2$ at the origin. Moreover the orders of the functions $z \mapsto f(0,0,z)$ and of $z \mapsto g(0,0,z)$ can be made arbitrarily large (in particular greater than~$2n$).
\end{prop}

Note also that, precisely as it happens in Sancho-Sanz example, for $p$ to be a persistent nilpotent singular point for the initial foliation, the latter must possess a formal separatrix $S = S_0$ through the point $p$ such that the singular point $p_n$ (selected by the transformed
separatrices $S_n$ in the sense that they correspond to the intersection of $S_n$ with the excetional divisor) is a nilpotent singular point for the corresponding foliation, for all $n \in \N$. Persistent nilpotent singular points also play a special role in the resolution theorem of \cite{MQ-P}. Namely, they appear as singularities associated
with a special type of $\Z /2\Z$-orbifold which, incidentally, require a weight~$2$ blow-up to be turned into elementary ones.
It is also worth pointing out that both statements are sharp in the sense that the well known example by Sancho and Sanz shows
the use of a weight~$2$ blow-up cannot be avoided (cf. Sections~$2$ and~$4$ of the same paper).

In particular, both Theorem~\ref{teo:A} and the resolution theorem, Theorem~$2$, in \cite{MQ-P} asserts the existence
of a birational model for $\fol$ where all singularities of $\fol$ are elementary except for finitely many ones that can be turned
into elementary singular points by means of a single blow-up of weight~$2$. In this sense, differences between these two theorems
are down to the way in which these rational models are constructed. Alternatively, Theorem~\ref{teo:A} can simply be regarded as a
new proof of the resolution theorem in \cite{MQ-P}.

In terms of the construction of the mentioned rational models, we briefly mention that McQuillan and Panazzolo work in the
category of {\it weighted blow-up}, along with the corresponding orbifolds, while in Theorem~\ref{teo:A} we restrict ourselves
as much as possible to the use of standard (i.e. unramified) blow-ups. Once again, additional information on these strategies
can be found in Section~2 of \cite{RR_Resolution}.

Throughout this section the term {\it blow-up}\, will
refer to {\it standard (i.e. homogeneous) blow-ups}. This applies, in particular, to the statement of Theorem~A. As to
blow-ups with weights, which are inevitably also involved in the discussion,
{\it these will be explicitly referred to as weighted (or ramified) blow-ups}.

Also, we will say that a (germ of) foliation $\fol$ can be
{\it resolved}\, if there is a sequence of blowing-ups as in~(\ref{blowup-resolution_globaldesingularization}) leading to
a foliation $\fol_n$ all of whose singularities are elementary. Similarly, a sequence of blowing-ups
as in~(\ref{blowup-resolution_globaldesingularization}) will be called a {\it resolution of $\fol$}\, if all the singular points
of $\fol_n$ are elementary. Whenever sequences of {\it weighted blow-ups}\, leading to a foliation having
only elementary singular points are considered, they may be referred to as a {\it weighted resolution of $\fol$}. With
this terminology, while every germ of foliation on $(\C^3,0)$ admits a weighted resolution, as follows from \cite{MQ-P}
or Theorem~A, the mentioned examples of Sancho and Sanz show that not all of them admit a resolution. Section~2
of \cite{RR_Resolution} contains a detailed discussion on the mutual interactions involving \cite{C-R-S}, \cite{MQ-P},
and our discussion revolving around Theorem~\ref{teo:A}.

\bigbreak

We can now go back towards our initial motivation, namely to germs of foliations $\fol$ on $(\C^3,0)$ that are
associated with a semicomplete vector field. Since
the notion of semicomplete singularity was introduced along with its first applications to the (global) study of complex
vector fields (\cite{Rebelo96}), it has been natural to ask whether all foliations in this class admit a resolution.
A special instance of this problem which is of interest in the
study of complex Lie group actions consists of asking whether the underlying foliation of a complete holomorphic
vector field (on some complex manifold of dimension~$3$) can be transformed into a foliation all of whose singular
points are elementary by means of a sequence of blow-ups as in~(\ref{blowup-resolution_globaldesingularization}).

To state our results concerning this special class of foliations, let us place ourselves once and for all in the context
of semicomplete vector fields. First, it is convenient to recall that a singularity of a holomorphic vector field $X$ is
said to be {\it semicomplete}\, if the integral curves of $X$ admit a maximal domain of definition in $\C$, cf.
Section~\ref{sec:sc_global_dynamics}. In particular, whenever $X$ is a {\it complete vector field}\, defined on a
complex manifold $M$, every singularity of $X$ is automatically semicomplete. The answer to the above question is
then provided by the following theorem:

\begin{teo}\cite{RR_Resolution}\label{teo:B}
	Let $X$ be a semicomplete vector field defined on a neighborhood of the origin in $\C^3$ and denote by $\fol$ the holomorphic
	foliation associated with $X$. Then one of the following holds:
	\begin{enumerate}
		\item The linear part of $X$ at the origin is nilpotent (non-zero).

		\item There exists a finite sequence of blow-ups maps along with transformed foliations
		$$
		\fol = \fol_0 \stackrel{\Pi_1}\longleftarrow \fol_1 \stackrel{\Pi_2}\longleftarrow \cdots
		\stackrel{\Pi_r}\longleftarrow \fol_r
		$$
		such that all of the singular points of $\fol_r$ are elementary. Moreover, each blow-up map $\Pi_i$ is centered in the
		singular set of the corresponding foliation $\fol_{i-1}$. In other words, the foliation $\fol$ can be resolved.
	\end{enumerate}
\end{teo}

Let us emphasize that item~(1) in Theorem~\ref{teo:B} means that the linear part of $X$ is (nilpotent) {\it non-zero from the outset}.
In other words, if the foliation $\fol$ associated with $X$ cannot be resolved, then $X$ has a non-zero nilpotent linear part
and this property is ``universal'' in the sense that it does not depend on any sequence of blow-ups/blow-downs carried out.
In particular, we can choose a ``minimal model'' for our manifold and the corresponding transform of $X$ will still have
non-zero nilpotent linear part at the corresponding point. Moreover, from Theorem~3 on \cite{RR_Resolution} about ``persistent
nilpotent singularities'', it is easy to obtain accurate normal forms for the vector field $X$.

Also, the statement of Theorem~\ref{teo:B} involves the linear part of the vector field $X$ rather than
the linear part of the associated foliation $\fol$. This makes for a stronger statement which is better
emphasized by Corollary~\ref{coro:C} below:

\begin{corollary}\cite{RR_Resolution}\label{coro:C}
	Let $X$ be a semicomplete vector field defined on a neighborhood of $(0,0,0) \in \C^3$
	and assume that the linear part of $X$ at the origin is equal to zero. Then item~(2) of Theorem~B holds.
\end{corollary}

More precisely, Theorem~\ref{teo:B} asserts that foliations associated with semicomplete vector fields in dimension~$3$ can be resolved
by a sequence of blow-ups centered in the singular set
except for a very specific case in which the vector field $X$ (and hence the foliation $\fol$) has a ``universal''
non-zero nilpotent linear part. As mentioned, these statements have the advantage of involving the vector field and not only the
underlying foliation. To clarify the meaning of this sentence, consider a holomorphic (semicomplete) vector field $X$
having the form $X = fY$, where $Y$ is another holomorphic vector field and $f$ is a holomorphic function. Whereas $X$ and
$Y$ induce the same singular foliation $\fol$, an immediate consequence
of Corollary~\ref{coro:C} is that $\fol$ must be as in item~(2) of Theorem~\ref{teo:B} {\it provided that $f$ vanishes at the origin}\,:
in fact, if $f$ and $Y$ are as indicated, then the linear part of $X$ vanishes at the origin at which $\fol$ is, indeed, singular
(clearly there is nothing to be proved if $\fol$ is regular). In other words, if $X = fY$ as above with $f(0,0,0) =0$ and
$X$ semicomplete, then the foliation associated with $X$ can certainly be resolved even if $Y$ has a nilpotent
singular point at the origin.

A few additional comments are needed to fully clarify the role of item~(1) in Theorem~B. First note that more accurate normal
forms are available for the vector fields in question: indeed, Theorem~3 of \cite{RR_Resolution} provides accurate normal forms
for all persistent nilpotent singular points. In addition, {\it not all}\, nilpotent vector fields giving rise to
persistent nilpotent singularities are semicomplete and, in this respect, the normal form provided by the mentioned
Theorem will further be refined.

Next, taking into account the global setting of complete vector fields, it is natural to wonder if there is, indeed,
{\it complete vector fields} inducing a foliation with singular points that cannot be resolved. As a consequence of
Theorem~\ref{teo:B}, such vector fields would definitely be pretty remarkable
since they must have a (non-zero) ``universal'' nilpotent singular point. To confirm that these global
situations do exist, however,  it suffices to note that the polynomial vector field
\[
Z = x^2 \partial /\partial x + xz \partial /\partial y + (y -xz) \partial /\partial z
\]
can be extended to a complete vector field defined on a suitable open manifold (details on Section~6 of \cite{RR_Resolution}).
As will be seen, the origin in the above coordinates constitutes a nilpotent singular point of $Z$ that cannot be resolved by
means of blow-ups as in item~(2) of Theorem~\ref{teo:B}, albeit this nilpotent singularity can be resolved by using a blow-up
centered at the (invariant) $x$-axis.

Finally, the question raised above about the existence of singularities as in item~(1) of Theorem~\ref{teo:B} in global
settings can also be asked in the far more restrictive case of holomorphic vector field defined on {\it compact manifolds}\,
of dimension~$3$. Owing to the compactness of the manifold, every such vector field is automatically complete. In this
setting, the methods used in the proof of Theorem~\ref{teo:B} easily yield:

\begin{corollary}\cite{RR_Resolution}\label{coro:D}
	Let $\fol$ be the foliation associated with a vector field $X$ globally defined
	on some compact manifold $M$ of dimension~$3$. Then every singular point of $\fol$ can be resolved.
\end{corollary}

Let us close this section with a couple of remarks inspired by some questions asked to us by A. Glutsyuk.
Essentially his questions concern resolution strategies with minimal number of (weighted) blow-ups which can also
be seen as an analogue of some questions previously considered in the context of Hironaka's theorem. In this respect,
it is clear that being able to work with weighted blow-ups, as opposed to standard ones, increases the chances of
reducing the number of blow-ups to resolve a given foliation. Indeed, it is easy to produce examples of this phenomenon
already in dimension~$2$ and in the context of Seidenberg's theorem. Hence, there is no chance that the strategy used in the proof
of Theorem~\ref{teo:A} will in general minimize the number of blow-ups required to resolve a given foliation. However,
we ignore if Panazzolo's algorithm \cite{P} has minimizing properties in the preceding sense.

A similar question directly motivated by the fact that in dimension~$2$ standard blow-ups suffice to resolve
any foliation, consists of trying to minimize the number of weighted
blow-ups needed to obtain the resolution. In this case, and at least for generic foliations, Theorem~\ref{teo:A} seems to provide
a satisfactory answer. Let us try to sketch an argument in this direction. As it follows from Theorem~3 of \cite{RR_Resolution},
persistent nilpotent singular points are naturally associated with certain formal separatrices (i.e. formal invariant curves)
having some special properties. Their ``position'' in the exceptional divisor obtained after finitely many blow-ups
is thus determined by the corresponding formal separatrices. In particular, it is possible to talk about these singularities being in
``general position'' for a given germ in an intrinsic way, i.e. independently of the use of any sequence of (standard) blow-ups. At
least when these singularities are in ``general position'' for a foliation $\fol$, then Theorem~\ref{teo:A} should minimize the number of
weighted blow-ups needed to turn $\fol$ into a foliation all of whose singular points are elementary. Indeed, each such singularity
requires at least one weighted blow-up to be turned into elementary singular points and each such blow-up can non-trivially
affect only one of these singularities thanks to the ``general position'' assumption. Thus the number of weighted blow ups
needed cannot be smaller than the number of persistent nilpotent singularities and the later is matched by the procedure in Theorem~A.
We ignore, however, if the ``general position assumption'' is really needed for this statement. Note that if there is
a foliation $\fol$ that can be resolved by using less weighted blow-ups than those prescribed in Theorem~A, then $\fol$ should
conceal at least two persistent nilpotent singularities so ``close'' to each other that they can both be turned into elementary
singular points by means of a same weighted blow-up.


	{\small\bibliography{cimart}}
 \EditInfo{ December 5, 2023}{February 10, 2024}{Ana Cristina Moreira Freitas, Carlos Florentino, Diogo Oliveira e Silva amd Ivan Kaygorodov
}
 
\end{document}